\documentclass[reqno,a4paper]{amsart}

\usepackage[utf8]{inputenc}
\usepackage[T1]{fontenc}
\usepackage[english]{babel}

\usepackage{latexsym}
\usepackage{amssymb}
\usepackage{dsfont}
\usepackage{mathtools}
\usepackage{enumitem}
\usepackage[hidelinks]{hyperref}
\usepackage{xcolor}

\allowdisplaybreaks[4]

\theoremstyle{plain}
\newtheorem{Theorem}{Theorem}[section]
\newtheorem*{Proposition*}{Proposition}
\newtheorem{Proposition}[Theorem]{Proposition}
\newtheorem{Lemma}[Theorem]{Lemma}

\newtheorem{Corollary}[Theorem]{Corollary}

\theoremstyle{definition}
\newtheorem{Remark}[Theorem]{Remark}

\numberwithin{equation}{section}

\renewcommand{\div}{\operatorname{div}}
\DeclareMathOperator{\curl}{curl}

\newcommand{\R}{\mathbb{R}}

\newcommand{\1}{\mathds{1}}

\newcommand{\cE}{\mathcal{E}}
\newcommand{\cD}{\mathcal{D}}
\newcommand{\cH}{\mathcal{H}}
\newcommand{\bH}{\mathbb{H}}

\newcommand{\abs}[1]{ \left\lvert#1\right\rvert} 
\newcommand{\norm}[1]{\left\lVert#1\right\rVert} 
\newcommand{\vertiii}[1]{{\left\vert\kern-0.25ex\left\vert\kern-0.25ex\left\vert #1
		\right\vert\kern-0.25ex\right\vert\kern-0.25ex\right\vert}}
\newcommand{\supp}{\operatorname{supp}}

\let\oldexists\exists \let\exists\relax \DeclareMathOperator{\exists}{\oldexists} 
\let\oldforall\forall \let\forall\relax \DeclareMathOperator{\forall}{\oldforall}

\makeatletter
\def\moverlay{\mathpalette\mov@rlay}
\def\mov@rlay#1#2{\leavevmode\vtop{%
		\baselineskip\z@skip \lineskiplimit-\maxdimen
		\ialign{\hfil$\m@th#1##$\hfil\cr#2\crcr}}}
\newcommand{\charfusion}[3][\mathord]{
	#1{\ifx#1\mathop\vphantom{#2}\fi
		\mathpalette\mov@rlay{#2\cr#3}
	}
	\ifx#1\mathop\expandafter\displaylimits\fi}
\makeatother
\newcommand{\cupdot}{\charfusion[\mathbin]{\cup}{\cdot}}

\newcommand\restr[2]{{
		\left.\kern-\nulldelimiterspace
		\vphantom{\big|}
		\right|_{#2}
}}

\newcommand{\tr}{\operatorname{tr}}

\newcommand*\dd{\mathop{}\!\mathrm{d}}

\newcommand{\e}{\mathrm{e}}
\newcommand{\pt}{\partial_t}

\newcommand{\inv}[1]{{#1}^{-1}}
\newcommand{\seq}{\hphantom{={}}}

\usepackage{multicol}

\renewcommand{\epsilon}{\varepsilon}
\newcommand{\Rsym}{\R^{3\times3}_{\mathrm{sym}}}

\newcommand{\+}{\hspace*{0.3pt}}
\newcommand{\ol}{\overline}
\newcommand{\Ra}{\mathrm{R}}
\newcommand{\Nu}{\mathrm{N}}

\begin{document}

\title[Interior Conductivity]{Exponential decay of the linear Maxwell system due to conductivity near the boundary}

\author[Richard Nutt and Roland Schnaubelt]{Richard Nutt$^1$ and Roland Schnaubelt$^2$}
\address{$^1$$^2$Karlsruhe Institute of Technology,
	Department of Mathematics,
	Englerstraße~2,
	76131 Karlsruhe,
	Germany}
\email{richard.nutt@kit.edu}
\email{schnaubelt@kit.edu}
\thanks{$^2$Corresponding author}

\thanks{Acknowledgment: Funded by the Deutsche Forschungsgemeinschaft (DFG, German Research Foundation) – Project-ID 258734477 – SFB 1173}



\begin{abstract}
	We study the anisotropic linear Maxwell system on a bounded domain $\Omega$ with perfectly conducting boundary conditions.
	It is damped via a 	conductivity $\sigma$ which is strictly positive on a collar at
	the boundary. We prove that solutions decay exponentially to 0,
	if the fields have no magnetic charges on $\Omega$ and
	no electric charges off the support of $\sigma$. Our approach relies on a splitting of the solution via a
	Helmholtz decomposition and an observability-type estimate for a related second-order system without
	charges, shown using Morawetz multipliers.
	Corresponding exact observability and controllability results are also established.
\end{abstract}

\maketitle

\section{Introduction}

The Maxwell equations are the fundamental laws of electro-magnetic theory. A non-zero conductivity $\sigma$
causes dissipation of energy and thus may lead to decay of solutions. For linear anisotropic materials we show
that the solutions converge exponentially to equilibria if $\sigma$ is strictly positive near the
boundary of the bounded domain $\Omega \subseteq \R^3$. It assumed that permittivity and permeability $\epsilon$ and $\mu$
satisfy the non-trapping condition \eqref{eq:non-trapping}. We allow for a multiply connected $\Omega$ and
disconnected $\partial\Omega$. We obtain exponential decay to 0 if there are no magnetic charges on $\Omega$ (which holds
in physics) and no electric charges on $\Omega\setminus \supp\sigma$. For the charge-free system with
$\sigma=0$, we also prove exact observability and controllability with respect to a collar at $\partial\Omega$.
In the case of anisotropic coefficients these seem to be the first results in this context if
the damping or observability region are not the full domain. Moreover, for exponential stability, so far only constant
$\epsilon$ and $\mu$ have been treated.

The Maxwell--Amp\`ere and Faraday equations relate the electric fields $E$ and $D$ with the magnetic ones $B$ and $H$ via
\[	\pt D(t,x) = \curl H(t,x) - J(t,x), \quad \pt B(t,x) = -\curl E(t,x),\qquad t\ge0, \; x\in\Omega,\]
where $J$ is a current density. One has to add constitutive relations that describe the reaction of the material to the fields.
We study anisotropic, instantaneous, time-independent, linear material laws
\[D(t,x)= \epsilon(x) E(t,x),\quad B(t,x)= \mu(x) H(t,x), \]
with \textit{permittivity} $\epsilon\in C^1\big(\ol{\Omega},\Rsym\big)$ and \textit{permeability}
$\mu \in C^1\big(\ol{\Omega},\Rsym\big)$ which are uniformly positive definite. In the absence of
exterior currents, Ohm's law yields
\[J(t,x) = \sigma(x)E(t,x) \]
for the non-negative \textit{conductivity} $\sigma \in L^\infty\big(\Omega,\Rsym\big)$.

We assume that there are no magnetic charges $\div(\mu H)$ and impose \textit{perfectly conducting boundary conditions},
thus excluding boundary conductivity. In this way we arrive at the Maxwell system
\begin{align}
	\begin{split}\label{eq:maxwell}
		\partial_t \big(\epsilon(x) E(t,x)\big) & = \curl H(t,x) - \sigma(x)E(t,x), \\
		\partial_t \big(\mu(x) H(t,x)\big) & = - \curl E(t,x),
	\end{split}
	 & & t\geq 0, \;x \in \Omega, \\
	\div(\mu(x)H(t,x)) & = 0, & & t\geq 0,\; x \in \Omega, \label{eq:maxwell2} \\
	E(t,x)\times\nu(x) & = 0, \quad \nu(x)\cdot\mu(x)H(t,x) = 0, & & t\geq 0, \; x\in\partial\Omega, \label{eq:maxwell3} \\
	E(0,x) & = E_0(x), \quad H(0,x)= H_0(x), & & x\in\Omega, \label{eq:maxwell4}
\end{align}
where $\nu$ is the outer unit normal of $\partial \Omega$. See \cite{ciarlet18}, \cite{cessenat1996},
or \cite{fabrizio2003} for a systematic treatment of the Maxwell equations from a mathematical point of view.

\smallskip

The anisotropy in the material laws is needed to describe non-trivial crystal structures, see e.g.\ \cite{BW99}.
Moreover, matrix-valued coefficients arise in the investigation of (even isotropic) quasilinear Maxwell systems, as in
\cite{lasiecka2019}, \cite{nutt2024} or \cite{pokojovy2020}. We stress that the anisotropic Maxwell system is
significantly stronger coupled than in the case of isotropic material laws with scalar coefficients.
For instance, if the coefficients are constant, the system can be reduced to decoupled scalar wave equations for
the components $E$ (if $\div(\epsilon E)=0$) or $H$ (if $\sigma =0$). In the isotropic case, the equations are coupled,
but only in lower order (if $\div(\epsilon E)=0$ or $\sigma =0$). Moreover, since anisotropic material laws
change the direction of the fields, the boundary conditions are much harder to handle in this setting.

The curl operator has a huge kernel consisting of gradient fields which prohibit regularity and compactness
properties. Divergence conditions may counteract the kernel. Physically they are encoded in the Gau{\ss}ian laws for
the charges. Besides $\div(\mu H)=0$, the Maxwell--Amp\`ere law in \eqref{eq:maxwell} yields
\begin{equation}\label{eq:gauss}
	\rho(t)\coloneqq\div(\epsilon E(t))= \div(\epsilon E_0) -\int_0^t \div(\sigma E(\tau))\dd \tau, \qquad t\ge0.
\end{equation}
Similarly one sees that the magnetic divergence and boundary conditions in \eqref{eq:maxwell2} and \eqref{eq:maxwell3}
are true if they are satisfied by $H_0$, cf.\ Theorem~5.2.5 in \cite{ciarlet18}. The time integral in \eqref{eq:gauss}
is a serious obstacle for the study of the long-time behavior. For isotropic coefficients one can reduce \eqref{eq:gauss}
to the scalar ODE $\partial_t\rho= - \frac{\sigma}{\epsilon}\rho- \nabla\frac{\sigma}{\epsilon}\cdot \epsilon E$
for fixed $x\in\Omega$. In the anisotropic case this approach fails.

Instead of interior damping, also boundary conductivities $\zeta$ have been studied. Here one sets $\sigma=0$
and replaces the boundary conditions \eqref{eq:maxwell3} by
\[ H\times \nu = \nu\times (\zeta E\times \nu),\]
see \cite{fabrizio2003}.
In this setting the electric charges $\rho$ in the domain vanish if $\rho(0)=0$. Moreover, this boundary condition can
improve the trace regularity compared to \eqref{eq:maxwell}--\eqref{eq:maxwell4}, see
\cite{schnaubelt2021} and the references therein. In this sense this case is simpler.

\smallskip

For the scalar wave equation there is a very rich and deep theory on decay properties caused by damping, whereas the
corresponding theory for the Maxwell system is far less developed. We refrain from giving
references for the wave case, but discuss a sample of the papers on the Maxwell equations.
For scalar coefficients exponential decay of solutions was shown in \cite{Ka94} for boundary damping with strictly
positive $\zeta$ under a non-trapping condition on a `substarlike' domain, see \cite{Ko94} for the case
of constant coefficients. Theorem~5.1 in \cite{phung2000} provides exponential decay for constant $\epsilon,\mu\in\R_+$ and
strictly positive scalar $\sigma$. Semilinear damping was treated in \cite{eller2002-2}. For $\sigma=0$ and scalar
coefficients, \cite{nicaise2005} establishes exact observability on a collar at $\partial \Omega$.
(See \cite{eller2002-2}, \cite{Ka94} and \cite{phung2000} for more results on observability and controllability.)
These theorems were proven by means of Morawetz multipliers and Helmholtz decompositions.
By \cite{phung2000}, for constant $\epsilon,\mu\in\R_+$ the solutions also decay exponentially
if $\sigma$ is strictly positive on its support and satisfies the geometric control condition from the wave case.
This proof is based on microlocal analysis.

For matrix-valued coefficients, so far exponential decay has been studied only assuming strict positivity of the
conductivity, see \cite{eller2019} and \cite{lasiecka2019} for interior damping as well as \cite{nutt2024} and
\cite{pokojovy2020} for the boundary case. The papers \cite{lasiecka2019}, \cite{nutt2024} and \cite{pokojovy2020} mainly deal
with nonlinear material laws. Earlier, in \cite{eller2007} boundary observability on $\partial\Omega$ was shown in the
anisotropic case. These papers also rely on Morawetz multipliers and Helmholtz decompositions. Weaker convergence properties
have been investigated, too, here we refer to \cite{NS25} for recent contributions and several references.

With the exception of Theorem~5.1 in \cite{phung2000}, in the cited works various topological constraints are
imposed, namely simple connectedness of $\Omega$ or connectedness of $\partial\Omega$, and partly even (variants of)
starshapedness of $\Omega$. The topological assumptions improve the mapping properties of the curl operator, as we recall
in Lemma~\ref{lem:curl-with-topology}. Starshapedness greatly simplifies the use of Morawetz multipliers for boundary
damping.

\smallskip

In our main result Theorem~\ref{thm:exp-decay} we show the uniform exponential decay to 0 of $L^2$-solutions
for matrix-valued $C^1$-coefficients assuming that $\sigma\in L^\infty$ is strictly positive on
$\omega\coloneqq\{\sigma>0\}$ containing a collar at the boundary. To exclude equilibria, besides the magnetic divergence
and boundary conditions we suppose that $\rho$ vanishes on $\Omega\setminus\omega$. The domain $\Omega$ is only required
to be multiply connected, in contrast to the previous literature. As in \cite{eller2007}, \cite{Ka94} or
\cite{nicaise2005} we assume the non-trapping condition \eqref{eq:non-trapping}
on $\epsilon$ and $\mu$, which says that the coefficients do not decay too rapidly in radial direction.
Heuristically, this property reduces backreflections, so that waves can reach the boundary in a sufficient way.
In Corollary~\ref{cor:main} we remove the condition on initial charges and the magnetic boundary condition.
In this framework we obtain uniform exponential convergence to the set of equilibria, which are related to non-trivial
charges. On the other hand, in the completely charge-free case without conductivity, Theorem~\ref{thm:obs} shows exact
observability of the Maxwell system by the electric fields on a collar at the $\partial \Omega$
and its exact controllability by charge-free current densities supported on the collar. As noted above, so far such
results with localized observations have been restricted to scalar coefficients, and for exponential stability
with localized damping even to constant $\epsilon$ and $\mu$.

In the next section we collect our assumptions, generation results and invertibility properties of the curl operator.
The main observability-type estimate only works in the charge-free case without conductivity, see system \eqref{eq:V_h}.
So we split $(E,H)$ into the solution $(V_h,W_h)$ of this system and other terms. To account for the charges on
$\omega$, we use a gradient field $\nabla p$ where $p\in H^1(\Omega)$ solves the elliptic PDE \eqref{def:p}. One can bound
$\partial_t \nabla p$ (but not $\nabla p$) in $L^2$ by the dissipation term $\|\sigma^{1/2} E\|_2$ which in turn can be
controlled via the energy equality. This approach goes back to \cite{phung2000} with $p\in H^1_0(\Omega)$ in the
constant coefficient case and for connected $\partial\Omega$ and $\supp \sigma$. Our system \eqref{def:p} takes care of
some the topological obstructions for the invertibility of $\curl$. To deal with the others, we have to restrict the
electric and magnetic fields in the Maxwell system to an invariant subsystem with finite codimension, see \eqref{def:X}.
The remaining part of the solution $(V_i,W_i)$ then solves the system \eqref{eq:V_nh} that incorporates the inhomogeneity
$-\sigma E-\epsilon \nabla\partial_t p$. Such a splitting was also used in \cite{nicaise2005} without $p$. To correct the
influence of $p(0)$ in these systems, one has to choose the initial values properly, see \eqref{def:wi0}. Finally, for
technical reasons some initial fields have to approximated by more regular ones.

In Proposition~\ref{prop:hom_obs} we show the crucial observability-type estimate for $V_h$. Dissipation terms
only occur for the electric part, so that it is useful to work with a second-order formulation involving only $V_h$.
Moreover, the behavior of $\nabla p$ and $\rho$ suggests to work with time derivatives.
So we first estimate $\partial_t V_h$ and pass only later to $V_h$ in Corollary~\ref{cor:obsV} by integration.
This result is then reformulated as Theorem~\ref{thm:obs} on observability and controllability. In a lengthy
calculation also involving vector analysis, Proposition~\ref{prop:hom_obs} is shown via a Morawetz multiplier and via
a tailor-made multiplier using Lax--Milgram (see Lemma~\ref{lem:Lax-Milgram}), as well as the energy equality for $V_h$.
For the reasoning, it is important that $\supp \sigma$ contains a collar at $\partial\Omega$ and that the coefficients
satisfy the non-trapping condition. For scalar coefficients similar arguments are found in \cite{nicaise2005}, but it is
quite sophisticated to extend them to the matrix-valued case.

In the last section we then estimate the other parts of $(E,H)$ in several steps. First, the energy of $\partial_t(E,H)$
is bounded by that of $(E,H)$ for which it is crucial that $\sigma$ is strictly positive on its support. This was done for
connected $\partial\Omega$ and $\omega$ in \cite{phung2000}, using the auxiliar function
$\nabla p$. In our more general setting we have to proceed differently, based on a lower estimate
for the curl operator in Lemma~\ref{lem:cE-E-estimate}.
One can control $\partial_t H$ in space-time by $\partial_t E$ and dissipation terms using the splitting of $E$
and energy-type estimates. The inhomogeneous part $(V_i,W_i)$ is handled just by means of Duhamel's formula, exploiting
properties of the initial data of \eqref{eq:V_nh}. Several of these inequalities on the time interval $[0,T]$ depend
on $T>0$. Nevertheless they can be combined into a proof of exponential stability employing a strategy that goes back
to \cite{Ha89} in the wave case. The convergence result for data with charges then follows by a projection argument.

\section{Notation, assumptions and auxiliary results}
\label{sec:1}

In this section we collect our main hypotheses and various basic properties. A bounded $C^2$-domain
$\Omega \subseteq \R^3$ is \emph{multiply connected} if there are disjoint $C^2$-surfaces $\Sigma_1,\dots, \Sigma_L$
such that $\Omega\setminus \bigcup_{l=1}^L \Sigma_l$ is simply connected, see \cite{ciarlet18}, \cite{cessenat1996}
or \cite{dautray1990}. Simple connected $\Omega$ are considered as the special case $L=0$.

\subsection*{Assumptions}
Throughout, we assume the following conditions.
\begin{enumerate}[leftmargin=21.77pt, label=\textup{(H)}]
	\item \label{hypothesis}
	 Let $\Omega\subseteq \R^3$ be open, bounded and multiply connected, and $\partial\Omega\in C^2$ have
	 connected components $\Gamma_0,\dots, \Gamma_K$.
	 Let $\epsilon, \mu \in C^{1}\big(\ol{\Omega}, \Rsym\big)$
	 be uniformly positive definite, $\sigma\in L^\infty\big(\Omega, \Rsym\big)$,
	 $\omega \coloneqq \{x\in \Omega\mid\sigma(x) >0\}$ have a $C^1$-boundary.
\end{enumerate}
The surface measure on $\partial\Omega$ or $\partial\omega$ is denoted by $\varsigma$.
We write $\upsilon \coloneqq \Omega \setminus \overline\omega$ and
$N_a \coloneqq \{ x \in \overline\Omega\mid\mathrm{dist}(x, \partial\Omega)<a\}$ for the collar of width $a>0$
at $\partial\Omega$. Our core assumption for exponential stability is that there are constants $a,\sigma_0>0$ such that
\begin{equation}\label{eq:conditions-on-conductivity}
	N_a\subseteq\omega \quad \text{and} \quad \sigma|_{\omega} \geq \sigma_0>0.
\end{equation}
In particular, this implies that $\partial \omega =\partial \upsilon \cupdot \partial \Omega$. $\upsilon \coloneqq \Omega \setminus \overline\omega$
Let $m(x)=x-x_0$ for some $x_0\in \R^3$. We also require the `non-trapping' condition
\begin{equation}\label{eq:non-trapping}
	\tilde\epsilon\coloneqq \epsilon+(m\cdot\nabla)\epsilon \geq \eta \epsilon, \qquad
	\tilde\mu \coloneqq \mu + (m\cdot \nabla)\mu \geq \eta \mu.
\end{equation}
for some $\eta > 0$; i.e, the coefficients cannot decay too strongly in radial direction.

We use the $L^2$-based Sobolev spaces $H^s$ on the domains and their boundaries.
The spaces $H(\curl)$ and $H(\div)$ are the maximal domains of the distributional operators $\curl$ and $\div$
in $L^2(\Omega)$ equipped with the graph norm, respectively. (We often omit the range spaces for vector fields etc.)
The kernels of these operators are written as $H(\curl0)$ and $H(\div0)$, and are endowed with the $L^2$-norm.
The tangential trace $\tr_t v= v|_{\partial \Omega}\times \nu$
and the normal trace $\tr_n v= v|_{\partial \Omega}\cdot\nu $
can be extended to maps from $H(\curl)$ and $H(\div)$ to $H^{-1/2}(\partial \Omega)$. Their kernels are denoted by
$H_{t0}(\curl)$ and $H_{n0}(\div)$, respectively, which are also the closure of test functions in the respective norms.
We also set $H_{t0}(\curl0)= H_{t0}(\curl)\cap H(\curl0)$ and $H_{n0}(\div0)=H_{n0}(\div) \cap H(\div0)$.
For $\alpha\in\{\epsilon,\mu\}$, the operator $\div_\alpha$ given by $\div_\alpha(v)=\div(\alpha v)$ with analogous
notations. Finally, we need the refinements
\begin{align*}
	H^\Sigma_{n0}(\div_\alpha 0) & \coloneqq \big\{u \in L^2(\Omega)^3 \bigm| \div_\alpha u =0, \, \tr_n(\alpha u) =0,\,
	\langle \tr_n (\alpha u), \1 \rangle_{H^{-\frac12}(\Sigma_l)}, \\ & \qquad \quad l\!\in\!\{ 1, \dots, L\}\big\}, \\
	H^\Gamma(\div_\alpha 0) & \coloneqq \big\{ u \in L^2(\Omega)^3 \bigm|\div_\alpha u =0, \,
	\langle \tr_n (\alpha u), \1 \rangle_{H^{-\frac12}(\Gamma_k)}, \, k \in\{ 0, \dots, K\}\big\}.
\end{align*}
For $\alpha= I$ we simply write $H^\Sigma_{n0}(\div 0)$ and $H^\Gamma(\div 0)$.
The divergence theorem yields $H^\Gamma(\div_\alpha 0)=H(\div_\alpha 0)$ if $\partial\Omega$ is connected,
and for simply connected $\Omega$ we have $H^\Sigma_{n0}(\div_\alpha 0)= H_{n0}(\div_\alpha 0)$.
See \cite{ciarlet18}, \cite{cessenat1996} or \cite{dautray1990} for these and related facts.

\subsection*{Generation results}
Let $L^2_\alpha(\Omega)$ be the space of measurable $f$ fulfilling $\alpha^{1/2} f \in L^2(\Omega)$, endowed
with the canonical norm. We define
\begin{equation} \label{def:X}
	X = \big\{ (e,h) \in L_{\epsilon}^2(\Omega) \times L^2_\mu(\Omega) \bigm|
	\epsilon E|_\upsilon \in H^\Gamma(\div 0,\upsilon),\;\mu h\in H^\Sigma_{n0}(\div 0)\big\},
\end{equation}
where $H^\Gamma(\div 0,\upsilon)$ is the variant of $H^\Gamma(\div 0)$ on $\upsilon$ etc. The operator
\begin{equation}\label{eq:generator}
	D(A) = \big(H_{t0}(\curl) \times H(\curl)\big) \cap X ,\qquad A = \begin{pmatrix}
		-\inv\epsilon \sigma & \inv\epsilon\curl \\ - \inv\mu\curl & 0
	\end{pmatrix},
\end{equation}
is well-defined by Theorem~6.1.4 in \cite{ciarlet18} and Proposition~IX.1.3 in \cite{dautray1990}, the latter
applied on $ \upsilon$. Note that the magnetic component of $(e,h)\in D(A)$ belongs to $H^1 (\Omega)$ by Proposition~6.1
of \cite{lasiecka2019}, and thus to $H^1 (\upsilon)$. One can then show its maximal dissipativity as in Lemma~2.1
of \cite{NS25}, so that $A$ generates a contractive $C_0$-semigroup $T(\cdot)$ on $X$.
Given $(E_0,H_0) \in D(A)$, we thus obtain a unique solution of \eqref{eq:maxwell}--\eqref{eq:maxwell4} satisfying
$\div (\epsilon E(t))=0$ on $\Omega \setminus \ol{\omega}$ for $t\ge0$ and
\begin{equation}\label{eq:solution}
	(E,H) \in C\big(\R_{\ge0}, D(A)\big) \cap C^1\big(\R_{\ge0}, X\big).
\end{equation}
We assume that $(E_0,H_0)$ belongs to $D(A)$ except for the proofs of Theorems~\ref{thm:obs} and \ref{thm:exp-decay}, where
we pass to general $(E_0,H_0) \in X$. Also in this case $(E(t),H(t))= T(t)(E_0,H_0)$ is called a solution of the system.

Conceptually, we exclude charges in the interior of $\Omega$ away from the support of the conductivity $\sigma$.
Furthermore, we assume that there are no charges on the connected components of the conductivity as in a capacitor,
since these generate static electric fields that do not decay in time. By Faraday's law we also exclude the
induction of currents in conductive loops along the boundaries $\partial\Sigma_i$ of the cuts of $\Omega$, for example
by an electric current passing through the `holes' of $\Omega$. (Think of a torus with an external current running
through the middle.)

The proof of Lemma~2.1 in \cite{NS25} implies that on $X_e=L_{\epsilon}^2(\Omega) \times L^2_\mu(\Omega)$
the extension $A_e$ of $A$ with domain $H_{t0}(\curl) \times H(\curl)$ generates a contraction semigroup $T_e(\cdot)$
which leaves $X$ invariant and coincides there with $T(\cdot)$. Hence, the divergence condition on
$\Omega \setminus \ol{\omega}$ for $E_0$, the divergence condition for $H_0$, and the boundary conditions on $\mu H_0$
are also invariant under the evolution for the system within $X_e$. Moreover, as we see in Lemma~\ref{lem:tildeA} the
kernel of $A_e$ is orthogonal to $X$ in $X_e$ so that our restrictions on charges and boundary conditions in
\eqref{def:X} exclude the stationary fields which obstruct exponential stability.

\subsection*{Energies}
Our arguments heavily use estimates of the energies'
\begin{equation}\label{def:ED} \begin{split}
		\cE(t) & \coloneqq \int_\Omega \left(\epsilon E(t) \cdot E(t) + \mu H(t) \cdot H(t)\right) \dd x
		= \norm{(E(t),H(t))}^2_X, \\
		\cD(t) & \coloneqq \int_\Omega \left(\epsilon \pt E(t) \cdot \pt E(t) + \mu \pt H(t) \cdot \pt H(t)\right) \dd x
		= \norm{(\pt E(t),\pt H(t))}^2_X,
	\end{split}\end{equation}
and in particular the following energy equality with dissipation related to $\sigma$.

\begin{Lemma}\label{lem:energy-decreasing}
	Let \ref{hypothesis} hold and $(E,H)$ as in \eqref{eq:solution} solve \eqref{eq:maxwell}. We then obtain
	\begin{align*}
		\cE(s) - \cE(t) & = 2\!\int_s^t\!\! \int_\Omega \big|\sigma^{1/2} E\big|^2 \dd x \dd \tau , \\
		\cD(s) - \cD(t) & = 2\!\int_s^t \!\!\int_\Omega \big|\sigma^{1/2} \pt E\big|^2 \dd x \dd \tau, \qquad t\ge s\ge0.
	\end{align*}
\end{Lemma}

\begin{proof}
	We only consider $\cD$ as the estimate for $\cE$ is shown in a similar, but simpler way. After regularizing $(E_0,H_0)$
	in $D(A^2)$, we pass to the time differentiated version of the Maxwell system \eqref{eq:maxwell}. Here the energy
	inequality is shown in a standard way.

	So, take fields $(E_0^n, H_0^n) \in D(A^2)$ converging to $(E_0, H_0)$ in $D(A)$ with respect to the graph norm.
	These initial values yield solutions
	\[ (E^n,H^n) \in C^2\big(\R_{\ge0}, X\big) \cap C^1\big(\R_{\ge0}, D(A)\big) \cap C\big(\R_{\ge0}, D(A^2)\big) \]
	of the equations
	\[ \pt^2 E^n = \inv \epsilon\curl \pt H^n - \inv\epsilon\sigma \pt E^n, \qquad \pt^2 H^n = -\inv\mu \curl\pt E^n.\]
	Then $(E^n, H^n)=T(\cdot)(E_0^n, H_0^n)$ tends to $(E,H)$ in the space $C_b\big(\R_{\ge0}, D(A)\big)$ of bounded functions,
	and $\pt (E^n, H^n)=T(\cdot)A(E_0^n,H_0^n)$ to $\pt(E,H)$ in $C_b\big(\R_{\ge0}, X\big)$.

	Let $\cD^n=\norm{(\pt E^n,\pt H^n)}^2_X$. The above system and integration by parts lead to
	\begin{align*}
		\pt \cD^n & = 2 \int_\Omega \big(\pt E^n \cdot\epsilon \pt^2E^n + \pt H^n \cdot\mu \pt^2 H^n\big) \dd x \\
		 & = 2 \int_\Omega \big(\pt E^n \cdot (\curl \pt H^n - \sigma \pt E^n) - \pt H^n \cdot\curl \pt E^n\big) \dd x \\
		 & = -2 \int_\Omega \big|\sigma^{1/2}\pt E^n\big|^2 \dd x,
	\end{align*}
	since $\tr_t\pt E^n = 0$. Integrating in time, we derive
	\[ \cD^n(s) - \cD^n(t) = 2\int_s^t\int_\Omega \big|\sigma^{1/2}\pt E^n\big|^2 \dd x \dd s.\]
	for $t\ge s\ge0$. The result follows in the limit $n \to \infty$.
\end{proof}

\subsection*{Properties of the curl operator}
We first list basic mapping properties of the gradient, divergence and curl. They can be found in Theorem~2.7 of
\cite{cessenat1996} or in \cite{dautray1990}, or follow similarly.

\begin{Lemma}\label{lem:mapping}
	In the following diagram the image of each operator is of finite codimension and contained in the kernel of the
	succeeding operator:
	\begin{align*}
		H^1(\Omega) & \xrightarrow{\nabla} \makebox[0pt][l]{$H(\curl)$}\hphantom{H_{t0}(\curl)} \xrightarrow{\curl}
		\makebox[0pt][l]{$H(\div)$}\hphantom{H_{n0}(\div)} \xrightarrow{\div} L^2(\Omega), \\
		H_c^1(\Omega) & \xrightarrow{\nabla} H_{t0}(\curl) \xrightarrow{\curl} H_{n0}(\div) \xrightarrow{\div} L^2(\Omega),
	\end{align*}
	where $H_c^1(\Omega)\coloneqq\big\{f\in H^1(\Omega)\bigm|\tr_{\Gamma_k}f\text{ constant for each } k\in\{0,\dots K\ \big\}
		\supset H^1_0(\Omega)$.
\end{Lemma}

The \emph{cohomology spaces}
\begin{align*}
	\bH_1 & \coloneqq \big\{ f \in L^2(\Omega) \bigm| \curl f = 0, \,\div f = 0, \,\tr_{n} = 0\big\}, \\
	\bH_2 & \coloneqq \big\{ f \in L^2(\Omega) \bigm| \curl f = 0, \, \div f = 0, \,\tr_{t} f = 0\big\}.
\end{align*}
contain functions which prevent the invertibility of $\curl$. In our setting they have finite dimensions, which are
determined by topological properties of the domain $\Omega$, see Remark~\ref{rem:top} and also Sections~3.2, 6.1 and 6.2
in \cite{ciarlet18}, Section~2.9 in \cite{cessenat1996}, or Section~IX.1.3 in \cite{dautray1990} for more details.

The next lemma collects results on mapping properties of $\curl$ from
Proposition~6.1.3, Theorem~6.1.4, Proposition~6.2.4 and Theorem~6.2.5 of \cite{ciarlet18}
and from Proposition~3, Proposition~4 and Remark~5 of Section~IX.1 in \cite{dautray1990} concerning the orthogonal
complements of the cohomology spaces.

\begin{Lemma}\label{lem:curl-with-topology}
	Let \ref{hypothesis} hold. We then have the surjectivity for
	\begin{enumerate}
		\item[a)] the `electric' curl
		 \begin{equation*}\begin{aligned} 
				 & \curl_E \colon H_{t0}(\curl) \cap H(\div_\epsilon 0) \to H^\Sigma_{n0}(\div0), \quad \text{with} \\
				 \Nu( & \curl_E) = \bH_2^\epsilon\coloneqq H_{t0}(\curl 0) \cap H(\div_\epsilon 0)
			 \end{aligned}\end{equation*}
		\item[b)] and the `magnetic' curl
		 \begin{equation*}\begin{aligned} 
				 & \curl_H \colon H(\curl) \cap H_{n 0}(\div_\mu 0) \to H^\Gamma(\div 0), \quad \text{with} \\
				 \Nu( & \curl_H) = \bH_1^\mu\coloneqq H(\curl0)\cap H_{n 0}(\div_\mu 0).
			 \end{aligned}\end{equation*}
	\end{enumerate}
	Furthermore, the image sets can be characterized by
	\begin{equation}\label{eq:char-bH}
		H^\Sigma_{n0}(\div0)= \bH_1^{\perp_{H_{n0}(\div 0)}}, \qquad H^\Gamma(\div 0)= \bH_2^{\perp_{H(\div 0)}}.
	\end{equation}
	In view of the next result we obtain the invertibility of
	\begin{alignat*}{3}
		\curl_E & \colon H_{t0}(\curl) \cap H^\Gamma(\div_\epsilon 0) & & \to H^\Sigma_{n0}(\div0) & & \qquad \text{in a) \ and} \\
		\curl_H & \colon H(\curl) \cap H_{n 0}^\Sigma(\div_\mu 0) & & \to H^\Gamma(\div 0) & & \qquad \text{in b)}.
	\end{alignat*}
\end{Lemma}

We have used the following description of the weighted cohomology spaces.

\begin{Lemma}\label{lem:cohomology-weighted}
	Let \ref{hypothesis} hold. For $u \in H(\div_\epsilon 0)$ the field $\epsilon u$ belongs to $\bH_2^{\perp_{L^2}}$
	if and only if $u$ is contained in $(\bH_2^\epsilon)^{\perp_{L^2_\epsilon}}.$
	Analogously, for $u \in H_{n0}(\div_\mu 0)$ the field $\mu u$ belongs to $\bH_1^{\perp_{L^2}}$ if and only if
	$u$ is contained in $(\bH_1^\mu)^{\perp_{L^2_\mu}}$.
\end{Lemma}

\begin{proof}
	For the first equivalence, we assume that $\epsilon u \in H^\Gamma(\div 0)$. Equation~(IX.1.61) in \cite{dautray1990}
	yields the decomposition
	\[ L^2(\Omega)^3 = \nabla H^1_0(\Omega) \oplus \bH_2 \oplus \curl H^1(\Omega). \]
	Since $H(\div 0)^{\perp_{L^2}}= \nabla H^1_0(\Omega)$ by (IX.1.42) in \cite{dautray1990},
	our assumptions and \eqref{eq:char-bH} imply that $\epsilon u$ is contained in $\curl H^1(\Omega)$; i.e.,
	$\epsilon u=\curl \Phi$ for some $\Phi\in H^1(\Omega)$. Hence, for all $h \in \bH_2^\epsilon$ we derive
	\[ (u,h)_{L^2_\epsilon} = (\epsilon u, h)_{L^2} = (\curl \Phi, h)_{L^2} = 0. \]
	The converse implication similarly follows from the decomposition
	\begin{equation}\label{proof:eq:weighted-L^2-decomposition}
		L^2_\epsilon(\Omega) = \nabla H^1_0(\Omega)\oplus_{L^2_\epsilon} \bH_2^\epsilon\oplus_{L^2_\epsilon}\epsilon^{-1}\curl H^1(\Omega),
	\end{equation}
	see Propositions~6.1.1 and 6.1.12 in \cite{ciarlet18} and Proposition~IX.1.3 in \cite{dautray1990}.

	The second assertion is shown analogously, using Propositions~6.2.1 and 6.2.12 in \cite{ciarlet18} and Proposition~IX.1.4
	in \cite{dautray1990}.
\end{proof}

\begin{Remark}\label{rem:top}
	Let \ref{hypothesis} hold. The dimension of the first cohomology space equals the `cutting number' $L$ of $\Omega$, i.e.,
	\[ \dim \bH_1 = \dim \Nu(\curl_H) = L. \]
	For simply connected $\Omega$ it is therefore trivial. Moreover, we have
	\[ \dim \bH_2 = \dim \Nu(\curl_E) = K.\]
	Thus, if $\partial\Omega$ is connected, the second cohomology space equals is trivial.
	See Proposition~2.8 in \cite{cessenat1996}, as well as Proposition~6.1.1 and Proposition~6.2.1 in \cite{ciarlet18}.
\end{Remark}

From now on we drop the subscript and simply write $\curl$ for both $\curl_E$ and $\curl_H$. We have to control the
normal trace of curls by tangential traces. To this aim we recall Lemma~4.8 of \cite{nutt2024}, see also Section~2.3
in \cite{cessenat1996}.
\begin{Lemma}\label{lem:surfaceCurl}
	For $f \in H^1(\Omega)$ we can estimate the normal trace of the curl by
	\[
		\norm{\nu \cdot\curl f}_{H^{-1}(\partial\Omega)} \lesssim \norm{\nu \times f}_{L^2(\partial\Omega)}.
	\]
\end{Lemma}

\subsection*{Helmholtz decomposition}
Our arguments rely on a splitting of the electric field
$E\in C^1\big(\R_{\ge0},L^2(\Omega)\big)\cap C\big(\R_{\ge0}, H_{t0}(\curl)\big)$ into a $\div_\epsilon$-free and
a $\curl$-free part. The following construction is inspired by Chapter~5 in \cite{phung2000}, where constant
$\epsilon$ and connected $\partial \Omega$ and $\omega$ were treated using $p\in H^1_0(\Omega)$.
We consider the elliptic problem
\begin{align}\label{def:p}
	 & \div(\epsilon \nabla p) = \div(\epsilon E) \quad \text{on } \Omega,\quad
	p = 0 \text{ on } \Gamma_0, \\
	 & p \text{ constant on } \Gamma_k\text{ \ and \ } \langle\tr_n\epsilon (\nabla p-E),\1\rangle_{H^{-\frac12}(\Gamma_k)}
	= 0, \ \ \forall\+ k\in \{1,\dots,K\}. \notag
\end{align}
For the weak formulation of the problem, we define the Hilbert space
\[\cH_c \coloneqq \big\{\varphi \in H^1(\Omega) \bigm| \varphi = 0 \text{ on }\Gamma_0, \,
	\varphi \text{ constant on } \Gamma_k, \ \forall\+ k\in\{1,\dots K\}\big\} \]
endowed with the norm of $H^1(\Omega)$.
We introduce the form $B[u,v] \coloneqq \int_\Omega \nabla u \cdot \epsilon \nabla v \dd x$ on $\cH_c$ and
set $\ell(v)=\ell_E(v) \coloneqq \int_{\Omega} \epsilon E \cdot \nabla v \dd x$. Then the weak formulation of \eqref{def:p}
reads as
\begin{equation}\label{eq:p-weak-formulation}
	B[p,\varphi] = \ell(\varphi) \quad \text{for all }\varphi \in \cH_c.
\end{equation}

Clearly, $B$ and $\ell$ are bounded on $\cH_c$, and $B$ is coercive by Theorem~13.6.9 in \cite{TW09}.
The Lax--Milgram lemma yields a unique solution $p=p(t) \in \cH_c$ of \eqref{eq:p-weak-formulation} with
\begin{equation}\label{est:p}
	\norm{p}_{H^1(\Omega)} \lesssim \norm{\ell}_{\cH_c^*} = \sup_{\norm{\varphi}_{\cH_c}=1} \abs{\ell(\varphi)}
	\lesssim \norm{\epsilon E}_{L^2(\Omega)}.
\end{equation}
We thus have bounded linear maps $L^2(\Omega)\to \cH_c^*$; $E\mapsto \ell_E$, and $L^2(\Omega)\to \cH_c$; $E\mapsto p$.
Note that $\epsilon E - \epsilon \nabla p$ belongs to $H(\div)$ with $\div(\epsilon E - \epsilon \nabla p)=0$,
since \eqref{eq:p-weak-formulation} holds in particular for $H^1_0(\Omega) \subseteq \cH_c$.
Hence, the last boundary condition in \eqref{def:p} is well-defined, and it follows by inserting
$\varphi_k \in \cH_c$ with $\varphi_k|_{\Gamma_j} = \delta_{kj}$ for $k\neq 0$ into \eqref{eq:p-weak-formulation} via
\[0= \int_\Omega (\epsilon E - \epsilon \nabla p)\cdot \nabla \varphi_k \dd x
	= \langle\tr_n (\epsilon E - \epsilon \nabla p),\1\rangle_{H^{-\frac12}(\Gamma_k)}. \]
Finally, the PDE in \eqref{def:p} is understood in $\cH_c^*$.

To pass to charge-free fields, we define
\begin{equation}\label{eq:def-V}
	V \coloneqq E -\nabla p\in C^1\big(\R_{\ge0},L^2(\Omega)\big)\cap C\big(\R_{\ge0}, H_{t0}(\curl)\big)
\end{equation}
using Lemma~\ref{lem:mapping} and \eqref{eq:solution}. By the properties of $p$, this function satisfies
\begin{equation}\label{eq:properties-V}
	\curl V = \curl E, \qquad	\epsilon V \in H^\Gamma(\div 0).
\end{equation}
For the last assertion, the condition $\langle\tr_n\epsilon V,\1\rangle_{H^{-\frac12}(\Gamma_k)}=0$
is clear by \eqref{def:p} for $k>0$. For $k=0$ it then follows from $\div(\epsilon V)=0$ which yields
$\langle\tr_n\epsilon V,\1\rangle_{H^{-\frac12}(\partial \Omega)}=0$.

In other main arguments we use the time derivatives of $p$ since they can be estimated solely by dissipation terms,
in contrast to $p$ in \eqref{est:p}. This fact corresponds to the behavior of the charges $\div(\epsilon E)$ in
\eqref{eq:gauss}, where the time derivative is better suited for the study of the long-term behavior, too.

\begin{Lemma}\label{lem:reg_of_ptp}
	Let $E$ be given by \eqref{eq:solution}. Then the solution $p$ of \eqref{def:p} belongs to	$C^2(\R_{\ge0},\cH_c)$ and satisfies
	\[\norm{\pt p}_{H^1} \lesssim \norm{\pt \nabla p}_{L^2} \lesssim \norm{\sigma E}_{L^2(\omega)}\quad \text{and} \quad
		\norm{\pt^2 p}_{H^1} \lesssim \norm{\pt^2 \nabla p}_{L^2}\lesssim \norm{\sigma \pt E}_{L^2(\omega)}.\]
	Moreover, $\epsilon \pt V$ is contained in $H^\Gamma(\div 0)$.
\end{Lemma}

\begin{proof}
	Because of \eqref{eq:solution}, we can continuously differentiate $t\mapsto\ell_{E(t)}$ in $\cH_c^*$ and thus $p$ in
	$\cH_c$. Equation \eqref{eq:maxwell} and Lemma~\ref{lem:curl-with-topology}\,b) also yield
	\begin{equation}\label{eq:ptp}
		\partial_t\ell_{E(t)}(\varphi) = \int_{\Omega} \epsilon \pt E \cdot \nabla \varphi \dd x
		= -\int_{\Omega} \sigma E \cdot \nabla \varphi \dd x
	\end{equation}
	for $\varphi\in \cH_c$. Estimate \eqref{est:p} then implies the first assertion, using also the Poincar\'e inequality in
	Theorem~13.6.9 in \cite{TW09}. The last claim now follows from \eqref{eq:properties-V}. Due to \eqref{eq:ptp},
	we can differentiate in time once more and derive the remaining part.
\end{proof}

\section{Observability estimate}
\label{sec:3}
Above we have split $E$ into the charge-free part $V$ and the curl-free part $\nabla p$. To obtain observability results,
we decompose $(V,H)$ further into fields $(V_h,W_h)$ and $(V_i,W_i)$ solving a homogeneous and an inhomogeneous system,
respectively. Here we combine and extend the approaches from \cite{nicaise2005} and \cite{phung2000}. The energy equality
in Lemma~\ref{lem:energy-decreasing} allows us to control dissipation terms, and thus the time derivatives of $\nabla p$
by Lemma~\ref{lem:reg_of_ptp}. Accordingly, we first estimate $\partial_t V_h$ and thus by-pass the problem that
the Gau\ss\ law \eqref{eq:gauss} does not provide uniform bounds in $t$.

In the stability analysis it will be crucial that $\pt V_{i}(0)$ and $\pt W_{i}(0)$ vanish, see
Lemma~\ref{lem:pt-V_nh-estimate}. To this aim, we choose suitable initial values for the inhomogeneous problem.
By Lemma~\ref{lem:reg_of_ptp} the derivative $\epsilon \pt V$ belongs to $H^\Gamma(\div 0)$, which implies
\begin{equation}\label{eq:inh}
	\sigma E + \epsilon \pt \nabla p = - \epsilon \pt E + \curl H + \epsilon \pt \nabla p
	= - \epsilon \pt V + \curl H \ \in H^\Gamma(\div 0)
\end{equation}
since $\curl H \in H^\Gamma(\div 0)$ due to Lemma~\ref{lem:curl-with-topology}\,b) as well as \eqref{eq:maxwell2}
and \eqref{eq:maxwell3}. Lemma~\ref{lem:curl-with-topology}\,b) then provides a field
$W_{i0}\in H(\curl)\cap H_{n 0}(\div_\mu 0) $ which is orthogonal to $\Nu(\curl_H)=\bH_1^\mu$ in $L^2_\mu$ and satisfies
\begin{equation}\label{def:wi0}
	\curl W_{i0} = \sigma E_0 + \epsilon \pt \nabla p(0).
\end{equation}

To derive the observability estimate for $\partial_t V_h$, we consider a second-order equation for our homogeneous problem
using the generator
\begin{equation}\label{eq:generator-hom}
	D(A_h) = \big(H_{t0}(\curl) \times H(\curl)\big) \cap X_h,\qquad A_h = \begin{pmatrix}
		0 & \inv\epsilon\curl \\ - \inv\mu\curl & 0
	\end{pmatrix}
\end{equation}
without conductivity acting in the charge-free subspace
\[ X_h \coloneqq \big\{ (v,w) \in L^2_\epsilon(\Omega) \times L^2_\mu(\Omega) \bigm|
	\epsilon v \in H^\Gamma(\div 0), \;\mu w \in H^\Sigma_{n0}(\div 0)\big\}.	\]
of $X$. The operator $A_h$ maps into $X_h$ by Lemma~\ref{lem:curl-with-topology}.
Its maximal dissipativity is a special case of that of $A$ in \eqref{eq:generator} for $\sigma=0$ and
$\omega =\emptyset$. By \eqref{def:X}, \eqref{eq:def-V}, \eqref{eq:properties-V} and
$W_{i0} \in H_{n 0}(\div_\mu 0)$, the fields $\big(E_0 - \nabla p(0), H_0-W_{i0}\big)$ belong to $D(A_h)$.

We can approximate $\big(E_0 - \nabla p(0), H_0-W_{i0}\big)$ in $D(A_h)$ by $(V_{h0}^{(n)}, W_{h0}^{(n)})\in D(A_h^2)$
by means of Yosida approximations. For these initial fields we can derive
a second-order problem. The homogeneous Maxwell system
\begin{align}\label{eq:V_h}
	\pt \big(\epsilon V^{(n)}_h\big) & = \curl W^{(n)}_h\+, & \hspace*{-1.5cm}
	\pt \big(\mu W^{(n)}_h\big) & = -\curl V^{(n)}_h\+,\notag \\
	\div \big(\epsilon V^{(n)}_h\big) & = 0, & \hspace*{-1.5cm}\div\big( \mu W^{(n)}_h\big) & = 0, \\
	\tr_t V^{(n)}_h & = 0, & \hspace*{-1.5cm}\tr_n(\mu W^{(n)}_h) & = 0, \notag \\
	V^{(n)}_h(0) & =V_{h0}^{(n)}, & \hspace*{-1.5cm} W^{(n)}_h(0) & =W_{h0}^{(n)}, \notag
\end{align}
then has the solution
\begin{equation}\label{eq:reg-VhWh}
	\big(V_h^{(n)},W_h^{(n)}\big) \in C^2\big(\R_{\ge0},X_h\big) \cap C^1\big(\R_{\ge0},D(A_h)\big)
	\cap C\big(\R_{\ge0},D(A_h^2)\big).
\end{equation}
It converges to the solution $(V_h, W_h)$ of \eqref{eq:V_h} with initial data $\big(E_0 - \nabla p(0), H_0-W_{i0}\big)$
in the space $C^1_b\big(\R_{\ge0},X_h\big) \cap C_b\big(\R_{\ge0}, D( A_h)\big)$. We stress that $\epsilon V_h^{(n)}$
is divergence-free. The system \eqref{eq:V_h} leads to the second-order equation
\begin{align}\label{eq:second-order-formulation}
	\pt \big(\epsilon \pt V_h^{(n)}\big) & = -\curl \big(\inv \mu \curl V_h^{(n)}\big), \notag \\
	\tr_t V^{(n)}_h & = 0, \quad \div (\epsilon V^{(n)}_h) =0, \\
	V_h^{(n)}(0) & =V^{(n)}_{h0}, \quad \pt V_h^{(n)}(0)=\inv \epsilon \curl W^{(n)}_{h0}\+.\notag
\end{align}
Corresponding to \eqref{eq:reg-VhWh}, in \eqref{eq:second-order-formulation} we assume that
\begin{equation}\label{eq:reg-Vh}\begin{split}
		V_h^{(n)} & \in C^2(\R_{\ge0},H^\Gamma(\div_\epsilon 0)) \cap C^1\big(\R_{\ge0},H_{t0}(\curl)\big), \\
		\mu^{-1}\curl V_h^{(n)} & \in C\big(\R_{\ge0},H(\curl)\big).
	\end{split}\end{equation}

To account for the missing terms in \eqref{eq:V_h}, we consider the inhomogeneous system
\begin{align}\label{eq:V_nh}
	\pt \big(\epsilon V_{i}\big) & = \curl W_{i} - \sigma E - \epsilon \pt \nabla p,\notag \\
	\pt \big(\mu W_{i}\big) & = -\curl V_{i}, \notag \\
	\div \big(\mu W_{i}\big) & = 0, \\
	\tr_t V_{i} & = 0, \quad \tr_n(\mu W_{i}) = 0, \notag \\
	V_{i}(0) & =0, \quad W_{i}(0)=W_{i0}\+. \notag
\end{align}
Since $\big({-}\epsilon^{-1}\sigma E- \pt\nabla p, 0\big)$ belongs to $C^1\big(\R_{\ge0},X_h\big)$ due to \eqref{eq:inh}
and $(0,W_{i0})$ to $D(A_h)$, this system has a unique solution
\begin{equation}\label{eq:Vi-reg}
	\big(V_{i}, W_{i}\big) \in C^1\big(\R_{\ge0}, X_h\big) \cap C\big(\R_{\ge0}, D( A_h)\big).
\end{equation}
For later use we record that the definition of $W_{i0}$ in \eqref{def:wi0} yields
\begin{equation}\label{eq:init-i}\begin{split}
		\pt (\epsilon V_{i})(0) & = \curl W_{i0} - \sigma E_0 - \epsilon \pt \nabla p(0) = 0, \\
		\pt (\mu W_{i})(0) & = -\curl V_{i}^{(n)}(0) = 0.
	\end{split}\end{equation}

Observe that the fields $\big(V_h + V_{i} + \nabla p , W_h + W_{i}\big)$
solve the original system \eqref{eq:maxwell}--\eqref{eq:maxwell4} with initial values
\[\big(E_0- \nabla p (0)+ \nabla p (0), H_0-W_{i0}+W_{i0}\big)=(E_0,H_0).\]
By uniqueness, the convergence noted after \eqref{eq:reg-VhWh} leads to the limit
\[ \big(V^{(n)}_h + V_{i} + \nabla p, W^{(n)}_h + W_{i}\big) \ \longrightarrow \
	\big(V_h + V_{i} + \nabla p, W_h + W_{i}\big) = (E,H)\]
in $C^1_b\big(\R_{\ge0}, X\big) \cap C_b\big(\R_{\ge0}, D(A)\big)$.

To simplify notation, we drop the superscripts in this section and write $\big(V_h^{(n)},W_h^{(n)}\big) = (V_h,W_h)$,
assuming that $(V_h, W_h)$ satisfies \eqref{eq:reg-VhWh}. The final result will then follow by the approximation argument
above, see Remark~\ref{rem:approx}. We next state the `energy' identity for $V_h$.

\begin{Lemma}\label{lem:energy}
	Let $V_h$ satisfy \eqref{eq:reg-Vh} and solve \eqref{eq:second-order-formulation}. We then obtain
	\begin{align*}
		\int_\Omega & \big(\epsilon \pt V_h(t) \cdot\pt V_h(t) + \inv\mu\curl V_h(t)\cdot\curl V_h(t)\big) \dd x \\
		 & = \int_\Omega \big(\epsilon \pt V_h(0) \cdot\pt V_h(0) + \inv\mu\curl V_h(0)\cdot\curl V_h(0)\big) \dd x, \qquad t\ge0.
	\end{align*}
\end{Lemma}

\begin{proof}
	System \eqref{eq:second-order-formulation} and integration by parts imply
	\begin{align*}
		0 & =\int_0^t \int_\Omega \big(\pt^2 (\epsilon V_h) + \curl \inv\mu\curl V_h\big) \cdot\pt V_h \dd x \dd t \\
		 & = \frac12\int_0^t \pt\int_\Omega \big(\pt(\epsilon V_h)\cdot\pt V_h + \inv\mu\curl V_h\cdot\curl V_h\big)\dd x\dd t \\
		 & = \frac12 \bigg[\int_\Omega \epsilon \pt V_h(t) \cdot\pt V_h(t) + \inv\mu\curl V_h(t)\cdot\curl V_h(t) \dd x \\
		 & \qquad-\int_\Omega \epsilon \pt V_h(0) \cdot\pt V_h(0) + \inv\mu\curl V_h(0)\cdot\curl V_h(0) \dd x\bigg]. \qedhere
	\end{align*}
\end{proof}

The following observability estimate for $\pt V_h$ is the core step in our arguments.

\begin{Proposition}\label{prop:hom_obs}
	Assume that \ref{hypothesis} and the non-trapping condition~\eqref{eq:non-trapping} hold. Let $V_h$ satisfy
	\eqref{eq:reg-Vh} and solve \eqref{eq:second-order-formulation}. Then for any $a > 0$ there exists a time $T_O >0$
	such that for $T\geq T_O$ we have
	\[ \int_\Omega \big(\abs{\pt V_h(0)}^2 + \abs{\curl V_h(0)}^2\big)\dd x
		\lesssim \int_0^T \!\!\int_{N_a} \abs{\pt V_h}^2 \dd x \dd t. \]
\end{Proposition}

For the proof of Proposition~\ref{prop:hom_obs} we need several auxiliary results. In the first one we use the Lax--Milgram
lemma to construct a multiplier, cf.\ \cite{nicaise2005}. We endow the space
$\cH\coloneqq H_{t0}(\curl) \cap H^\Gamma(\div_\epsilon 0)$ with the norm $\norm{\cdot}_{H(\curl)}$.
By Lemma~\ref{lem:curl-with-topology} the operator
$\curl_E \colon \cH \to H^\Sigma_{n0}(\div 0)$ is boundedly invertible.

\begin{Lemma}\label{lem:Lax-Milgram}
	Take a map $\tilde \vartheta \in L^\infty(\Omega)$ with $\tilde \vartheta = 1$ on $N_{a/2}$ and
	$\supp \tilde\vartheta \subseteq N_a$ for some $a>0$, and $f \in L^2(\Omega)$. Then there is a unique solution
	$w \in \cH \subseteq H^1(\Omega)$ of
	\begin{equation}\label{eq:Lax-Milgram1}
		\forall \psi \in \cH:\qquad \int_\Omega \inv\mu \curl w \cdot\curl \psi \dd x = \int_\Omega \tilde \vartheta f \psi \dd x.
	\end{equation}
	It satisfies
	\begin{equation}\label{eq:Lax-Milgram2}
		\|w\|_{H^1}\lesssim\norm{w}_{\cH} \lesssim	\|\tilde\vartheta f\|_{L^2}\+.
	\end{equation}
\end{Lemma}

\begin{proof}
	Theorem~A.6\,a) of \cite{nutt2024} (with $v=0$ and $u \in H_{t0}(\curl)$) yields
	\begin{equation*}
		\norm{u}_{H^1} \lesssim \norm{u}_{L^2} + \norm{\curl u}_{L^2} + \norm{\div (\epsilon u)}_{L^2}\+.
	\end{equation*}
	Since the divergence of $\epsilon u$ vanishes for $u \in \cH$ and the curl is invertible, we see that
	\begin{equation*}
		\norm{u}_{H^1} \lesssim \norm{\curl u}_{L^2} \qquad \text{for } u \in \cH\+.
	\end{equation*}
	Hence, the bilinear form on the left-hand side is coercive and bounded.
	The assertion then follows from the Lax--Milgram lemma, using that $L^2(\Omega) \hookrightarrow \cH^\star$.
\end{proof}

For a time dependent $f \in C^1\big(\R_{\ge0},L^2\big)$ we obtain a solution $w(t)$ of \eqref{eq:Lax-Milgram1} for
each $t\ge0$. The regularity of $f$ transfers to $w$ as shown next.

\begin{Remark}\label{rem:regularity-w}
	Let $f \in C^1\big(\R_{\ge0},L^2(\Omega)\big)$ in Lemma~\ref{lem:Lax-Milgram}. We then obtain
	\begin{enumerate}
		\item[a)] $w\in C^1\big(\R_{\ge0},\cH\big)$ and
		\item[b)] $\norm{w(t)}^2_{H^1} \lesssim \int \tilde\vartheta \abs{f(t)}^2\dd x$, \
		 $\norm{\pt w(t)}^2_{H^1} \lesssim \int \tilde\vartheta \abs{\pt f(t)}^2\dd x$.
	\end{enumerate}
\end{Remark}

\begin{proof}
	We first check that $w$ is $C^1$. Its Lipschitz continuity follows from \eqref{eq:Lax-Milgram2} and linearity via
	\[ \big\|\tfrac1h \big(w(t+h)-w(t) \big)\big\|_{H^1}\lesssim \big\|\tfrac1h \big(f(t+h) - f(t) \big)\big\|_{L^2}\+.\]
	Hence, $w$ is differentiable in $H^1(\Omega)$ for a.e.\ $t\ge0$. We can now differentiate \eqref{eq:Lax-Milgram1} in $t$
	a.e., obtaining
	\[ \int_\Omega \inv\mu \curl \pt w \cdot\curl \psi \dd x = \int_\Omega \tilde \vartheta \pt f \psi \dd x. \]
	As above we then infer the continuity of $\pt w$ from \eqref{eq:Lax-Milgram2}, which also yields the second estimate.
	The first one follows from \eqref{eq:Lax-Milgram2} for $w$.
\end{proof}

We now turn to the proof of the core result of this section.

\begin{proof}[Proof of Proposition~\ref{prop:hom_obs}]
	\begin{enumerate}[label = \arabic*), wide, labelwidth=0pt]
		\item Let $T>0$, $\chi \in C^1(\ol{\Omega})$, and $m(x) = x-x_0$ for some $x_0 \in \R^3$. We multiply
		 \eqref{eq:second-order-formulation} with the Morawetz multiplier $\chi m \times \curl V_h$ obtaining
		 \begin{align}\label{eq:mora0}
			 0 & =\int_0^T\!\!\int_\Omega \big(\pt(\epsilon \pt V_h) + \curl\inv\mu\curl V_h\big)\cdot\big(\chi m\times\curl V_h\big)\dd x \dd t\notag \\
			 & = \int_0^T \!\!\int_\Omega \!\pt\big(( \epsilon \pt V_h)\cdot (\chi m \!\times\! \curl V_h)\big) \dd x \dd t
			 - \int_0^T \!\!\int_\Omega (\epsilon \pt V_h) \cdot\big(\chi m \!\times\! \curl \pt V_h\big) \dd x \dd t \notag \\
			 & \seq +\int_0^T\!\!\int_\Omega \inv\mu \curl V_h \cdot\curl \big(\chi m \times \curl V_h) \dd x \dd t \notag \\
			 & \seq -\int_0^T\!\!\oint_{\partial \Omega} \inv\mu \curl V_h \cdot\big( \nu \times (\chi m \times \curl V_h)\big)\dd\varsigma \dd t\notag \\
			 & \eqqcolon I_i +I_t + I_c + I_c^{\partial}\+.
		 \end{align}
		 In the above four integrands we have to eliminate the second derivatives of $V_h$, which is easy for $I_i$.
		 In the following we freely use standard formulas from vector analysis. Also integrating by parts
		 and exploiting $\tr_t \pt V_h = 0$, we compute
		 \[ I_t =\int_0^T \!\!\int_\Omega \curl \pt V_h\cdot(\chi m \times \epsilon \pt V_h) \dd x \dd t
			 = \int_0^T\!\! \int_\Omega \pt V_h\cdot\curl(\chi m \times \epsilon \pt V_h) \dd x \dd t.\]
		 Since $\div (\epsilon V_h) =0$, the integrand can be rewritten as
		 \begin{align*}
			 \curl & (\chi m \times \epsilon \pt V_h) \\
			 & = \big(\epsilon \pt V_h \cdot\nabla\big)(\chi m) - \chi (m \cdot \nabla)(\epsilon \pt V_h)
			 + \chi m \pt\div (\epsilon V_h) - \epsilon \pt V_h \div(\chi m) \\
			 & = \big(\epsilon \pt V_h\cdot \nabla \chi\big)m + \chi \epsilon \pt V_h - \chi (m \cdot \nabla)(\epsilon \pt V_h)
			 - (\epsilon \pt V_h)\nabla \chi \cdot m - 3\chi \epsilon \pt V_h.
		 \end{align*}
		 Only the middle term involves second derivatives of $V_h$. In this summand of the integrand we can take out the gradient via
		 \[ \chi (m\cdot \nabla)(\epsilon \pt V_h )\cdot \pt V_h
			 =\tfrac12\chi(m\cdot \nabla)(\epsilon \pt V_h \cdot \pt V_h)
			 +\tfrac12 \chi\big((m\cdot\nabla)\epsilon\big) \pt V_h \cdot \pt V_h\]
		 by means of the product rule. The remaining problematic part of $I_t$ is now integrated by parts resulting in
		 \begin{align*}
			 {-} \int_0^T\!\! \int_\Omega & \frac12\chi (m\cdot \nabla)(\epsilon \pt V_h \cdot \pt V_h)\dd x \dd t \\
			 & = \int_0^T \!\! \int_\Omega\Big( \frac12(\nabla \chi \cdot m) (\epsilon \pt V_h \cdot \pt V_h)
			 + \frac32 \chi(\epsilon \pt V_h \cdot \pt V_h) \Big) \dd x \dd t \\
			 & \quad - \int_0^T\!\! \oint_{\partial \Omega} \frac12\chi \nu \cdot m (\epsilon \pt V_h \cdot \pt V_h) \dd \varsigma \dd t.
		 \end{align*}
		 Recalling $\tilde \epsilon = \epsilon + (m\cdot \nabla) \epsilon$ from \eqref{eq:non-trapping}, we conclude
		 \begin{align}\label{eq:It}
			 I_t & = \int_0^T\!\!\int_\Omega\!\Big[{-}\frac12\chi\tilde\epsilon \pt V_h \!\cdot\!\pt V_h
				 - \frac12 (\nabla\chi \!\cdot\! m)(\epsilon\pt V_h\!\cdot\!\pt V_h)
			 + \pt V_h \cdot (\epsilon \pt V_h \!\cdot\! \nabla \chi)m \Big]\mathrm{d}x \dd t \notag \\
			 & \seq - \int_0^T \!\!\oint_{\partial \Omega} \frac12 \chi \nu \cdot m (\epsilon \pt V_h\! \cdot\! \pt V_h) \dd \varsigma\dd t.
		 \end{align}

		 We next treat $I_c$ in a similar way, reformulating it as
		 \begin{align*}
			 I_c & = \int_0^T\!\! \int_\Omega \inv \mu \curl V_h\cdot \big[\!\div \curl (V_h) \chi m - \div (\chi m)\curl V_h) \\
			 & \qquad \qquad + (\curl V_h \cdot \nabla)(\chi m) - (\chi m \cdot \nabla) \curl V_h\big] \dd x \dd t.
		 \end{align*}
		 The last term is the only problematic one. As above it is equal to
		 \begin{align*}
			 {-} & \int_0^T \!\!\int_\Omega (\inv \mu \curl V_h)\cdot \big((\chi m \cdot \nabla) \curl V_h\big)\dd x \dd t \\
			 & = \frac12\int_0^T \!\!\int_\Omega\Big[(\chi m \cdot \nabla\inv\mu ) \curl V_h \cdot\curl V_h
				 -\chi (m\cdot \nabla) (\inv\mu\curl V_h \cdot \curl V_h) \Big]\textrm{d}x \dd t.
		 \end{align*}
		 We integrate by parts the last term and obtain
		 \begin{align} \label{eq:Ic}
			 I_c & =\int_0^T \!\!\int_\Omega \inv\mu \curl V_h \cdot \Big((\curl V_h \cdot \nabla \chi)m + \chi \curl V_h
			 - \frac12 \div(\chi m) \curl V_h\Big) \mathrm{d} x \dd t\notag \\
			 & \quad+ \frac12 \int_0^T\!\! \int_\Omega \big((\chi m \cdot \nabla)\inv\mu\big)\curl V_h \cdot\curl V_h\dd x \dd t\notag \\
			 & \quad -\frac12 \int_0^T\!\! \oint_{\partial\Omega}\chi (m\cdot\nu)(\inv\mu\curl V_h \cdot\curl V_h) \dd \varsigma \dd t \notag \\
			 & = \int_0^T \!\!\int_\Omega \Big[{-}\frac12 \chi {\tilde \mu}_{-1} \curl V_h \cdot \curl V_h
			 + \inv \mu \curl V_h \cdot \big((\curl V_h \cdot \nabla \chi)m\big) \notag \\
			 & \qquad \quad - \frac12 (m \cdot \nabla \chi) (\curl V_h \cdot \inv\mu \curl V_h)\Big]\textrm{d}x \dd t\notag\notag \\
			 & \quad - \frac12 \int_0^T \!\!\oint_{\partial\Omega}\chi (m\cdot\nu)(\inv\mu\curl V_h \cdot\curl V_h) \dd \varsigma \dd t,
		 \end{align}
		 Note that \eqref{eq:non-trapping} implies ${\tilde \mu}_{-1} \coloneqq \inv \mu - m \cdot \nabla \inv \mu \geq \eta \inv \mu$.

		 Using Lemma~\ref{lem:surfaceCurl}, the last integral in \eqref{eq:mora0} can be rewritten as
		 \begin{align}\label{eq:Ip}
			 I_c^{\partial} & = -\int_0^T\!\!\oint_{\partial\Omega}\inv\mu\curl V_h\cdot\big(\nu\times(\chi m\times\curl V_h)\big)\dd\varsigma\dd t\notag \\
			 & = \int_0^T \!\!\oint_{\partial \Omega} \nu \cdot \big(\inv\mu \curl V_h \times(\chi m \times \curl V_h)\big) \dd \varsigma \dd t\notag \\
			 & = \int_0^T \!\!\oint _{\partial \Omega} \chi(\nu \cdot m) (\inv\mu \curl V_h \cdot \curl V_h) \dd \varsigma \dd t.
		 \end{align}

		 Equations \eqref{eq:mora0}, \eqref{eq:It}, \eqref{eq:Ic} and \eqref{eq:Ip} then yield the core identity
		 \begin{align}\label{eq:mora}
			 \int_0^T\!\! & \int_\Omega\frac12\chi\big(\tilde \epsilon \pt V_h \cdot \pt V_h + \inv{\tilde\mu}\curl V_h \cdot\curl V_h\big)\dd x\dd t \\
			 & =\left[\int_\Omega \epsilon \pt V_h \cdot (\chi m \times \curl V_h) \dd x\right]_0^T\notag \\
			 & \quad -\int_0^T \!\!\int_\Omega\frac12 (\nabla \chi \cdot m) \big(\epsilon \pt V_h \cdot \pt V_h
			 + \inv \mu \curl V_h \cdot \curl V_h\big)\dd x\dd t \notag \\
			 & \quad + \int_0^T\!\! \int_\Omega \big((m \cdot \pt V_h)(\pt \epsilon V_h \cdot \nabla \chi)
			 + (m \cdot \inv \mu \curl V_h)(\curl V_h \cdot \nabla \chi)\big)\dd x\dd t\notag \\
			 & \quad + \int_0^T \!\!\oint_{\partial \Omega} \frac12 \chi (m \cdot \nu)
			 \big( \inv\mu \curl V_h \cdot \curl V_h-\epsilon \pt V_h \cdot \pt V_h\big) \dd \varsigma \dd t.\notag
		 \end{align}

		\item We first take $\chi=1$. Then \eqref{eq:mora} and the non-trapping condition \eqref{eq:non-trapping} imply
		 \begin{align*}
			 \eta \int_0^T \!\! & \int_\Omega \big(\abs{\pt V_h}^2 + \abs{\curl V_h}^2\big) \dd x \dd t \\
			 & \lesssim \int_\Omega \abs{\pt V_h(T)}\abs{\curl V_h(T)} \dd x + \int_\Omega \abs{\pt V_h(0)}\abs{\curl V_h(0)} \dd x \\
			 & \qquad+\int_0^T\!\!\oint_{\partial\Omega}\big(\nu\cdot m)( \inv\mu\curl V_h\cdot\curl V_h- \epsilon\pt V_h\cdot\pt V_h\big)\dd\varsigma\dd t.
		 \end{align*}
		 Next, for $\chi\in C^1(\ol{\Omega})$ with $\chi=1$ on $\partial\Omega$ and support in $N_{a/4}$, Equation~\eqref{eq:mora} yields
		 \begin{align*}
			 \Big|\int_0^T\!\!
			 & \oint_{\partial\Omega}(\nu\cdot m)\big(\inv\mu\curl V_h\cdot\curl V_h-\epsilon\pt V_h\cdot\pt V_h\big)\dd\varsigma\dd t\Big| \\
			 & \lesssim \int_\Omega \abs{\pt V_h(T)}\abs{\curl V_h(T)}\dd x+ \int_\Omega \abs{\pt V_h(0)}\abs{\curl V_h(0)} \dd x \\
			 & \seq + \int_0^T\!\!\int_{N_{a/4}} \!\!\big(\abs{\pt V_h}^2 + \abs{\curl V_h}^2\big) \dd x \dd t.
		 \end{align*}
		 The last two estimates lead to
		 \begin{align}\label{eq:proof1}
			 \int_0^T\!\! \int_\Omega \big(|\pt V_h|^2 + |\curl V_h|^2\big) \dd x \dd t
			 \lesssim \int_0^T\!\!\int_{N_{a/4}} \big(\abs{\pt V_h}^2 + \abs{\curl V_h}^2 \big)\dd x \dd t\notag \\
			 +\int_\Omega \abs{\pt V_h(T)} \abs{\curl V_h(T)} \dd x+ \int_\Omega \abs{\pt V_h(0)} \abs{\curl V_h(0)} \dd x.
		 \end{align}

		\item\label{proof:obs-3} In light of Lemma~\ref{lem:energy} we only have to estimate the curl term on $N_{a/4}$ to obtain the desired
		 result, see \eqref{eq:proof2}. To this aim, take $\vartheta \in C^1(\ol{\Omega})$ with $\supp \vartheta \subseteq N_{a/2}$ and
		 $ \vartheta = 1$ on $N_{a/4}$. Equation~\eqref{eq:second-order-formulation} and integration by parts yield
		 \begin{align*}
			 0 & = \int_0^T\!\! \int_\Omega (\pt^2( \epsilon V_h) + \curl (\inv\mu \curl V_h))\cdot(\vartheta V_h) \dd x \dd t \\
			 & = -\int_0^T \!\!\int_\Omega \epsilon\pt V_h \cdot \pt( \vartheta V_h)\dd x \dd t
			 + \left[\int_\Omega ( \epsilon\pt V_h) \cdot \vartheta V_h \dd x\right]_0^T \\
			 & \seq+ \int_0^T \!\!\int_\Omega\inv \mu \curl V_h \cdot \big(\vartheta \curl V_h+ \nabla\vartheta \times V_h\big)\dd x\dd t \\
			 & \seq + \int_0^T\!\!\oint_{\partial \Omega} (\nu \times \inv\mu\curl V_h)\cdot(\vartheta V_h)\dd \varsigma \dd t.
		 \end{align*}
		 The boundary integral vanishes by \eqref{eq:second-order-formulation} since
		 \[(\nu \times \inv\mu\curl V_h)\cdot(\vartheta V_h) = (\vartheta V_h \times \nu)\cdot(\inv\mu\curl V_h) = 0\]
		 on $\partial \Omega$. 	As $\mu^{-1}\ge \eta/\|\mu\|_\infty$, for any $\delta>0$ we derive
		 \begin{align*}
			 & \int_0^T \!\!\int_{N_{a/4}} \abs{\curl V_h}^2 \dd x \dd t \leq \int_0^T \!\!\int_{N_{a/2}} \vartheta \abs{\curl V_h}^2 \dd x \dd t \\
			 & \lesssim \delta \int_0^T\!\!\int_{N_{a/2}} \abs{\curl V_h}^2\dd x \dd t + c_\delta \int_0^T\!\!\int_{N_{a/2}} \abs{V_h}^2 \dd x\dd t
			 + \int_0^T \!\! \int_{N_{a/2}} \abs{\pt V_h}^2 \dd x \dd t \\
			 & \quad + \int_{\Omega} \big(\abs{\pt V_h(T)}^2 + \abs{V_h(T)}^2\big) \dd x
			 + \int_{\Omega} \big(\abs{\pt V_h(0)}^2 + \abs{V_h(0)}^2\big) \dd x.
		 \end{align*}
		 We insert this inequality into \eqref{eq:proof1} and absorb the curl term fixing a small $\delta>0$. Using also
		 Lemma~\ref{lem:energy}, it follows
		 \begin{align}\label{eq:proof3}
			 \int_0^T & \!\! \int_\Omega \big(\abs{\pt V_h}^2 + \abs{\curl V_h}^2\big) \dd x \dd t \\
			 & \lesssim \int_0^T \!\!\int_{N_{a/2}} \big(\abs{V_h}^2 + \abs{\pt V_h}^2 \big)\dd x \dd t
			 + \int_{\Omega} \big(\abs{\pt V_h(0)}^2 + \abs{\curl V_h(0)}^2 \big)\dd x.\notag
		 \end{align}

		\item We have to get rid of the new term involving $V_h$. Taking $w\in H_{t0}(\curl)$ from Lemma~\ref{lem:Lax-Milgram}
		 with $f=V_h$ as a multiplier, like in step \ref{proof:obs-3} we derive
		 \begin{align*}
			 0 & =\int_0^T \!\!\int_\Omega \big(\pt^2 (\epsilon V_h) + \curl(\inv\mu \curl V_h)\big)\cdot w \dd x \dd t \\
			 & = \int_0^T \!\! \int_\Omega \big((\inv\mu \curl V_h)\cdot \curl w - \epsilon\pt V_h \cdot\pt w\big)\dd x \dd t
			 + \left[\int_\Omega \epsilon\pt V_h \cdot w \dd x\right]_0^T.
		 \end{align*}
		 As $V_h$ belongs to $\cH$ by \eqref{eq:reg-Vh}, Equation \eqref{eq:Lax-Milgram1}
		 with $V_h=\psi$ then yields
		 \begin{align*}
			 0 & = \int_0^T\!\! \int_\Omega \big(\tilde \vartheta \abs{V_h}^2 - \epsilon\pt V_h \cdot\pt w \big)\dd x \dd t
			 + \left[\int_\Omega \epsilon \pt V_h \cdot w \dd x\right]_0^T\,.
		 \end{align*}
		 The properties of $\tilde\vartheta$, Remark~\ref{rem:regularity-w}, Lemma~\ref{lem:energy} and
		 Lemma~\ref{lem:curl-with-topology}\,a) then imply
		 \begin{align*}
			 \int_0^T\!\! \int_{N_{a/2}} \abs{V_h}^2 \dd x \dd t
			 & \lesssim \delta \int_0^T \!\!\int_\Omega \abs{\pt V_h}^2 \dd x \dd t
			 + c_\delta \int_0^T\!\! \int_\Omega \tilde\vartheta \abs{\pt V_h}^2 \dd x \dd t \\
			 & \seq +\! \int_\Omega\!\! \big(|\pt V_h(T)|^2 + |V_h(T)|^2\big) \mathrm{d} x
			 + \!\int_\Omega\!\!\big(|\pt V_h(0)|^2 + |V_h(0)|^2\big) \mathrm{d} x \\
			 & \lesssim\delta\int_0^T \!\!\int_\Omega\abs{\pt V_h}^2\dd x \dd t
			 + c_\delta\int_0^T \!\!\int_\Omega \tilde\vartheta \abs{\pt V_h}^2 \dd x \dd t \\
			 & \seq+ \int_\Omega \big(\abs{\pt V_h(0)}^2 + \abs{\curl V_h(0)}^2\big) \dd x.
		 \end{align*}
		 We insert this inequality in \eqref{eq:proof3}. For a fixed small $\delta>0$, the space-time term without localization
		 can be absorbed by the left-hand side resulting in
		 \begin{align}\label{eq:proof4}
			 \int_0^T \!\!\int_\Omega & \big(|\pt V_h|^2 + |\curl V_h|^2 \big)\dd x \dd t \\
			 & \lesssim \int_0^T \!\!\int_{N_a} \abs{\pt V_h}^2 \dd x \dd t
			 + \int_\Omega \big(|\pt V_h(0)|^2 + |\curl V_h(0)|^2 \big)\dd x.\notag
		 \end{align}

		\item In a final step we simplify the time integral on the left-hand side of \eqref{eq:proof4}
		 by means of Lemma~\ref{lem:energy}, and obtain
		 \begin{align}\label{eq:proof2}
			 T\int_\Omega \big(|\pt V_h(0)|^2 + |\curl V_h(0)|^2 \big)\dd x
			 & \leq c_0\int_\Omega \big(|\pt V_h(0)|^2 + |\curl V_h(0)|^2 \big)\dd x \notag \\
			 & \seq + c_1\int_0^T \!\!\int_{N_a} \abs{\pt V_h}^2 \dd t \dd x
		 \end{align}
		 for some constants $c_j>0$. Taking $T>c_0$, we conclude
		 \begin{equation*}
			 \int_\Omega \big(|\pt V_h(0)|^2 + |\curl V_h(0)|^2\big) \dd x
			 \leq \frac{c_1}{T-c_0}\int_0^T \!\!\int_{N_a} \abs{\pt V_h}^2 \dd t \dd x. \qedhere
		 \end{equation*}
	\end{enumerate}
\end{proof}

One can also show such an observability estimate for $V_h$ instead of $\pt V_h$; that is, for the solutions to a
homogeneous Maxwell system ($\sigma = 0$) with divergence free initial values. As in \cite{nicaise2005} we pass to
an antiderivative in the proof.

\begin{Corollary}\label{cor:obsV}
	Assume that \ref{hypothesis} and the non-trapping condition~\eqref{eq:non-trapping} hold. Let $(V_h,W_h)$ satisfy
	\eqref{eq:reg-Vh} and solve \eqref{eq:V_h}.
	For each $a > 0$ there exists a time $T_O >0$ such that for $T\geq T_O$ we have
	\[\int_\Omega \big(|V_h(0)|^2 + |W_h(0)|^2\big) \dd x \lesssim \int_0^T \int _{N_a} \abs{V_h}^2 \dd x \dd t.\]
\end{Corollary}

\begin{proof}
	Since $\mu W_h(0) \in H^\Sigma_{n0}(\div0)$, Lemma~\ref{lem:curl-with-topology}\,a) yields
	a field $u_0$ in $H_{t0}(\curl) \cap H^\Gamma(\div_\epsilon0)$ with $\curl u_0 = - \mu W_h(0).$ 	We set
	\[ u(t) \coloneqq u_0 + \int_0^t V_h(s) \dd s, \qquad t\ge0,\]
	obtaining $u(0) =u_0$ and $\pt u(0) = V_h(0)$. Note that $u$ satisfies \eqref{eq:reg-Vh}.
	We check that $u$ solves \eqref{eq:second-order-formulation}
	with different initial values. The definition of $u_0$ and \eqref{eq:V_h} yield
	\begin{align*}
		\pt^2 (\epsilon u(t)) & = \pt(\epsilon V_h(t) = \curl W_h(t)
		= \curl\Big(\inv \mu \Big(\mu W_h(0) + \int_0^t \pt (\mu W_h(s))\dd s \Big)\Big) \\
		 & = -\curl \Big(\inv \mu \curl\Big( u_0 + \int_0^t V_h(s) \dd s \Big)\Big)= -\curl \big(\inv \mu \curl u(t)\big), \\
		\tr_t u(t) & = \tr_tu_0 + \int_0^t \tr_t V_h(s) \dd s = 0, \\
		\div(\epsilon u(t)) & = \div(\epsilon u_0) + \int_0^t \div(\epsilon V_h(s)) \dd s = 0.
	\end{align*}
	The corollary now follows from Proposition~\ref{prop:hom_obs} as $\pt u=V_h$ and $\curl u= -\mu W_h$.
\end{proof}

\begin{Remark}\label{rem:approx}
	In Lemma~\ref{lem:energy}, Proposition~\ref{prop:hom_obs} and Corollary~\ref{cor:obsV} we have assumed
	that $(V_h(0), W_h(0)) = \big(V_{h0}^{(n)}, W_{h0}^{(n)}\big)$ belongs to $D(A_h^2)$ and so the solutions
	$(V_h, W_h) = \big(V_h^{(n)},W_h^{(n)}\big)$ satisfy \eqref{eq:reg-VhWh}. By the approximation argument discussed after
	\eqref{eq:V_nh}, the lemma and the proposition can be extended to $(V_h(0), W_h(0))$ in $D(A_h)$. Similarly, one
	derives the corollary	for $(V_h(0), W_h(0)) \in X_h$.
\end{Remark}

We reformulate the above corollary as (exact) observability and controllability of \eqref{eq:maxwell}--\eqref{eq:maxwell4}
with $\div(\epsilon E)=0$ and $\sigma=0$. In this charge-free case we have $p=0$ and $(V_h,W_h)=(E,H)$,
see \eqref{def:p}, \eqref{def:wi0}, \eqref{eq:V_h} and \eqref{eq:V_nh}. The equivalence of observability and controllability
is shown in Theorem~11.2.1 in \cite{TW09}, noting that $A_h$ is skew-adjoint in
$X_h=H^\Gamma(\div_\epsilon 0)\times H_{n0}^\Sigma(\div_\mu 0)$, see \eqref{eq:generator-hom}.

\begin{Theorem}\label{thm:obs}
	Assume that \ref{hypothesis} with $\sigma=0$ and the non-trapping condition~\eqref{eq:non-trapping} hold.
	Let $a>0$. Then there is a time $T_O>0$ and a constant $c>0$ such that for $(E_0,H_0)\in X_h$ we have
	\[\int_\Omega \big(|E_0|^2 + |H_0|^2\big) \dd x \le c \int_0^{T_O} \int _{N_a} \abs{E}^2 \dd x \dd t,\]
	where $(E,H)\in C(\R_{\geq 0} ,X_h)$ solves \eqref{eq:maxwell}--\eqref{eq:maxwell4} with $\sigma =0$.

	Moreover, for each $T\ge T_O$ and $(E_1,H_1)\in X_h$ we find a current $J\in L^2((0,T)\times N_a)$ with values
	in $H^\Gamma(\div_\epsilon 0)$ such that
	$(\hat E(T), \hat H(T))=(E_1,H_1)$ for the solution $(\hat E, \hat H)\in C(\R_{\geq 0} ,X_h)$
	of \eqref{eq:maxwell}--\eqref{eq:maxwell4} with $\sigma E$ replaced by $J$.
\end{Theorem}

\section{Exponential decay}

We now present the second main result of the paper.

\begin{Theorem}\label{thm:exp-decay}
	Let \ref{hypothesis}, \eqref{eq:conditions-on-conductivity} and \eqref{eq:non-trapping} hold.
	Then there exist constants $M \ge 1$ and $\omega>0$ such that
	\[ \|(E(t),H(t))\|_{L^2} \leq M\e^{-\omega t} \|(E_0,H_0)\|_{L^2}, \qquad t\geq 0,\]
	for the solution $(E,H)$ of \eqref{eq:maxwell}--\eqref{eq:maxwell4} with initial value $(E_0,H_0) \in X$.
\end{Theorem}

Our proof relies on Lemma~\ref{lem:cE-estimate} which estimates the usual energy through
time derivatives. This is needed to control the inhomogeneous part $(V_i,W_i)$ via Lemma~\ref{lem:pt-V_nh-estimate}.
Theorem~5.4 of \cite{phung2000} provides such an inequality for constant $\epsilon$ and $\mu$ and
for connected $\partial\Omega$ and $\omega$. The argument there is based on the splitting $E=V+\nabla p$.
We could not extend this approach to our setting and proceed in a more direct way using the following estimate
of $E$ on $\upsilon$.

\begin{Lemma}\label{lem:cE-E-estimate}
	For $E \in H^\Gamma(\div_\epsilon 0, \upsilon) \cap H(\curl,\upsilon)$, we have
	\[ \norm{E}_{L^2(\upsilon)} \lesssim \norm{\curl E}_{L^2(\upsilon)} + \norm{\tr_t E}_{H^{-1/2}(\partial \upsilon)}.\]
\end{Lemma}

\begin{proof}
	By Proposition~IX.1.3 of \cite{dautray1990} the field $E$ belongs to
	$\epsilon^{-1}\curl H^1(\upsilon)= H^\Gamma(\div_\epsilon 0, \upsilon)$, see also Lemma~\ref{lem:curl-with-topology}.
	In view the 	decomposition~\eqref{proof:eq:weighted-L^2-decomposition}, we can thus
	compute its norm in $L^2_\epsilon(\upsilon)$ by testing with $\varphi=\epsilon^{-1}\curl\Phi$ for $\Phi\in H^1(\upsilon)$.
	Theorem~3.4.1 in \cite{ciarlet18} allows us to choose the vector potential $\Phi$ such that
	$\norm{\Phi}_{H^1(\upsilon)} \lesssim \norm{\epsilon\varphi}_{L^2(\upsilon)}$. It follows
	\begin{equation}\label{proof:L^2-norm}
		\norm{E}_{L^2_\epsilon(\upsilon)}
		= \sup_{\substack{\varphi\in H^\Gamma(\div_\epsilon 0, \upsilon),\\ \norm{\varphi}_{L^2_\epsilon(\upsilon)} =1}}		\int_\upsilon E \cdot \epsilon\varphi \dd x
		\lesssim \sup_{\substack{\Phi\in H^1(\upsilon),\\ \norm{\Phi}_{H^1(\upsilon)} =1}}
		\int_\upsilon E \cdot \curl \Phi \dd x.
	\end{equation}
	Integration by parts then yields
	\begin{align*}
		\int_\upsilon E \cdot \curl \Phi \dd x & = \int_{\upsilon} \curl E \cdot \Phi \dd x
		+ \langle \tr_t E, \Phi \rangle_{H^{-1/2}(\partial \upsilon)} \\
		 & \leq \norm{\curl E}_{L^2(\upsilon)} \norm{\Phi}_{L^2(\upsilon)}
		+ \norm{\tr_t E}_{H^{-1/2}(\partial \upsilon)} \norm{\Phi}_{H^1(\upsilon)}
	\end{align*}
	We infer the assertion by inserting this inequality into \eqref{proof:L^2-norm}.
\end{proof}

The next proof relies on the assumption that $\sigma$ either vanishes
or is uniformly positive.

\begin{Lemma}\label{lem:cE-estimate}
	Let \ref{hypothesis} and \eqref{eq:conditions-on-conductivity} hold. Take a solution $(E,H)$ of
	\eqref{eq:maxwell}--\eqref{eq:maxwell4} as in \eqref{eq:solution}. The energies $\cE$ and $\cD$ from \eqref{def:ED}
	then satisfy
	\begin{equation*}
		\cE(T) \lesssim \cD(T), \qquad T\ge0.
	\end{equation*}
\end{Lemma}

\begin{proof}
	We adopt some ideas from the proof of Theorem~5.4 in \cite{phung2000}. Condition \ref{hypothesis},
	the Maxwell equations~\eqref{eq:maxwell}, \eqref{eq:maxwell3} and integration by parts lead to
	\begin{align}\label{proof:eq-omega}
		\sigma_0 \int_{\omega} \abs{E(t)}^2 \dd x & \leq \int_\Omega \sigma E(t) \cdot E(t)
		= \int_\Omega \big({-}\epsilon\pt E(t) +\curl H(t)\big)\cdot E(t) \dd x \notag \\
		 & = - \int_\Omega \epsilon \pt E(t)\cdot E(t) - \int_\Omega H(t)\cdot\mu \pt H(t)\dd x \notag \\
		 & \lesssim \sqrt{\cE(t)}\sqrt{\cD(t)}.
	\end{align}
	Using also \eqref{eq:solution}, Lemma~\ref{lem:curl-with-topology}\,b) and \eqref{eq:maxwell}, we then estimate
	the magnetic field by
	\begin{equation}\label{est:H}
		\norm{H(t)}^2_{L^2}\lesssim \norm{\curl H(t)}^2_{L^2}
		\lesssim \norm{\epsilon \pt E(t)}^2_{L^2} + \norm{\sigma E(t)}^2_{L^2} \lesssim \cD(t) + \sqrt{\cE(t)\cD(t)}.
	\end{equation}
	To control $E$ on $\upsilon$, we recall
	\[ \norm{E(t)}_{L^2(\upsilon)} \lesssim \norm{\curl E(t)}_{L^2(\upsilon)}^2
		+ \norm{\tr_t E(t)}_{H^{-1/2}(\partial \upsilon)}. \]
	from Lemma~\ref{lem:cE-E-estimate}. Since $E \in H(\curl)$, Proposition~2.2.32 in \cite{ciarlet18}
	shows that $\tr_{t,\partial \upsilon} E = \tr_{t,\partial \omega} E$ on $\partial \upsilon$.
	From the usual trace estimate and Equation~\eqref{eq:maxwell} we thus deduce
	\begin{align*}
		\norm{E(t)}_{L^2(\upsilon)}
		 & \!\lesssim \norm{\curl E(t)}_{L^2(\upsilon)} \!+ \norm{\tr_t E(t)}_{H^{-\frac12}(\partial \omega)}
		\!\lesssim \norm{\mu \pt H(t)}_{L^2} \!+ \norm{E(t)}_{H(\curl,\omega)} \\
		 & \!\lesssim \norm{\mu \pt H(t)}_{L^2} + \norm{E(t)}_{L^2(\omega)}.
	\end{align*}
	Combined with \eqref{proof:eq-omega} and \eqref{est:H}, we arrive at
	\[
		\cE (t) \le c_1\sqrt{\cE(t)}\sqrt{\cD(t)} + c_2\cD(t)\le \tfrac12 \cE(t) + C \cD(t)
	\]
	which yields the assertion.
\end{proof}

The time derivative of magnetic field can be estimated by the electric one using our Helmholtz decomposition.
We proceed similar as in the proof of Lemma~5.2 of \cite{phung2000}, see also Proposition~4.4 in \cite{lasiecka2019}.
Here and below the estimates depend on the end time, which fortunately does not cause problems in the main argument.
\begin{Lemma}\label{lem:pt-H-estimate}
	Let \ref{hypothesis} hold. Take a solution $(E,H)$ of \eqref{eq:maxwell}--\eqref{eq:maxwell4} as in \eqref{eq:solution}.
	For $T >0$ we then obtain
	\[ \int_0^T \!\!\int_\Omega \abs{\pt H}^2 \dd x \dd t \lesssim
		(1+T) \int_0^T \!\!\int_\Omega \abs{\pt E}^2 \dd x \dd t + \int_0^T \int_\Omega \sigma E \cdot E\dd x \dd t. \]
\end{Lemma}

\begin{proof}
	Fix $\Phi \in C_c^\infty((0,T))$ such that $0\leq \Phi \leq 1$ and $\Phi = 1$ on $[\tfrac13T,\tfrac23T]$. In \eqref{eq:def-V}
	and \eqref{eq:properties-V}, we have decomposed $E = V+\nabla p$ with $V \in H_{t0}(\curl) \cap H^\Gamma(\div_\epsilon0)$.
	Furthermore (after regularization) $E$ solves the second-order problem
	\[\epsilon \pt^2 E = -\curl(\inv \mu \curl E) - \sigma \pt E. \]
	Starting from \eqref{eq:maxwell} and integrating by parts, we then compute
	\begin{align*}
		\int_0^T\!\! \int_\Omega \Phi^2 \mu \pt H \cdot \pt H \dd x \dd t\,
		 & = \int_0^T \!\!\int_\Omega \Phi^2 \curl E \cdot\inv \mu \curl E \dd x \dd t \\
		 & = \int_0^T \!\!\int_\Omega \Phi^2 \curl V \cdot \inv \mu \curl E \dd x \dd t \\
		 & = \int_0^T \!\!\int_\Omega \Phi^2 V \cdot \curl \big(\inv \mu \curl E\big) \dd x \dd t \\
		 & = -\int_0^T \!\!\int_\Omega\Phi^2 V\cdot \big(\epsilon \pt^2 E + \sigma \pt E\big)\dd x \dd t \\
		 & = \int_0^T \!\!\int_\Omega \big((2\Phi \pt\Phi )V + \Phi^2 \pt V\big )\cdot\epsilon \pt E \dd x \dd t \\
		 & \seq - \int_0^T \!\!\int_\Omega\Phi^2 V\cdot \sigma \pt E\dd x \dd t.
	\end{align*}
	Lemma~\ref{lem:curl-with-topology}\,a) yields $\norm{V}_{L^2} \lesssim \norm{\curl V}_{L^2}$. We also insert
	$\partial_t V= \pt E-\pt\nabla p$.	H\"older's inequality thus implies
	\begin{multline*}
		\int_0^T \!\!\int_\Omega \Phi^2 \mu \pt H \cdot \pt H \dd x \dd t
		\lesssim\int_0^T\Big[\norm{\Phi'}_\infty\big(\delta\big\|\Phi\mu^{-1/2}\curl V\big\|_{L^2}^2 + \frac1\delta\norm{\pt E}_{L^2}^2\big)\\
			+\big(\norm{\pt \nabla p}_{L^2}^2 + \norm{\pt E}_{L^2}^2\big) + \big(\delta \big\|\Phi \mu^{-1/2}\curl V\big\|_{L^2}^2
			+ \frac1\delta \norm{\sigma \pt E}_{L^2}^2\big)\Big] \dd t.
	\end{multline*}
	As $\curl V = -\mu \pt H$ by \eqref{eq:maxwell}, we can absorb the curl terms by the left-hand side.
	Lemma~\ref{lem:reg_of_ptp} now leads to
	\begin{align*}
		\int_0^T \!\!\int_\Omega \Phi^2 \mu \pt H \cdot \pt H \dd x \dd t
		 & \lesssim \int_0^T \int_\Omega \abs{\pt E}^2 \dd x \dd t + \int_0^T \int_\Omega \abs{\sigma E}^2 \dd x \dd t.
	\end{align*}
	By means of the energy estimates from Lemma~\ref{lem:energy-decreasing}, we conclude
	\begin{align*}
		\int_0^T \!\!\int_\Omega \mu \pt H \cdot \pt H \dd x \dd t
		 & \leq T\cD(0) = T\left(\cD(T) + 2 \int_0^T \sigma \pt E \cdot \pt E\right) \dd x \dd t \\
		 & \leq 3 \tfrac T3 \,\cD\left(\tfrac23 T\right) + 2 T \int_0^T \sigma \pt E \cdot \pt E \dd x \dd t \\
		 & \lesssim \int_{T/3}^{2T/3} \abs{\mu^\frac12\pt H}^2 \dd x \dd t + (1+T) \int_0^T \abs{\pt E}^2 \dd x \dd t \\
		 & \lesssim (1 + T)\int_0^T \int_\Omega \abs{\pt E}^2 \dd x \dd t + \int_0^T \int_\Omega \sigma E\cdot E \dd x \dd t.\qedhere
	\end{align*}
\end{proof}

We next treat the inhomogeneous part of the fields. Duhamel's formula and our choice of the initial values lead to
the following estimate.

\begin{Lemma}\label{lem:pt-V_nh-estimate}
	Let \ref{hypothesis} hold and $(V_{i},W_{i})$ as in \eqref{eq:Vi-reg} solve \eqref{eq:V_nh}. We then obtain
	\begin{align*}
		\int_0^T\!\!\int_\Omega \big(\abs{\pt V_{i}}^2 + \abs{\pt W_{i}}^2 \big)\dd x \dd t
		 & \leq CT^2 \int_0^T\int_\Omega \abs{\sigma \pt E}^2 \dd x \dd t\,.
	\end{align*}
\end{Lemma}

\begin{proof}
	We use the generator $A_h$ of $T_h(\cdot)$ from \eqref{eq:generator-hom}. Equation~\eqref{eq:V_nh} and the subsequent
	comments imply that
	\[ \begin{pmatrix} V_{i}(t)\\ W_{i}(t) \end{pmatrix} = T_h(t) \begin{pmatrix} 0\\W_{i0} \end{pmatrix}
		- \int_0^t T_h(s) \begin{pmatrix}\inv\epsilon\sigma E(t-s) + \pt \nabla p(t-s) \\ 0\end{pmatrix}\dd s.\]
	Because of \eqref{eq:solution} and Lemma~\ref{lem:reg_of_ptp} we can differentiate this formula in $L^2$ with respect to $t$,
	where the resulting initial values vanish due to \eqref{eq:init-i}; i.e.,
	\[ \pt \begin{pmatrix} V_{i}(t) \\ W_{i}(t) \end{pmatrix}
		= -\int_0^t T_h(t-s) \begin{pmatrix}\inv\epsilon\sigma \pt E(s) + \pt^2 \nabla p(s) \\ 0\end{pmatrix}\dd s. \]
	Lemma~\ref{lem:reg_of_ptp} and H\"older's inequality now yield
	\begin{align*}
		\int_0^T \norm{\begin{pmatrix}\pt V_{i} \\ \pt W_{i} \end{pmatrix}}_{L^2}^2 \dd t
		 & \lesssim \int_0^T \left( \int_0^t \norm{\sigma \pt E}_{L^2} \dd s \right)^2 \dd t
		\lesssim \int_0^T t \int_0^t \norm{\sigma \pt E}^2_{L^2} \dd s \dd t \\
		 & \leq CT^2 \int_0^T\!\!\int_\Omega \abs{\sigma \pt E}^2 \dd x \dd t.\qedhere
	\end{align*}
\end{proof}

The above estimates lead to the core inequality.
\begin{Proposition} \label{prop:decay}
	Assume that \ref{hypothesis}, \eqref{eq:conditions-on-conductivity} and \eqref{eq:non-trapping} hold. Let
	$(E,H)$ as in \eqref{eq:solution} solve \eqref{eq:maxwell}--\eqref{eq:maxwell4}, and $T_O>0$ be given by
	Proposition~\ref{prop:hom_obs}. Then there exists a constant $\gamma\in[0,1)$ such that
	\[ \cE(T) + \cD(T) \leq \gamma \big(\cE(0) + \cD(0)\big), \qquad T\geq \max\{T_O,1\}. \]
\end{Proposition}

\begin{proof}
	Lemmas~\ref{lem:cE-estimate} and \ref{lem:energy-decreasing} imply
	\[ T(\cE(T) + \cD(T)) \lesssim T \cD(T) \lesssim \int_0^T\int_\Omega \big(\abs{\pt E}^2 + \abs{\pt H}^2\big)\dd x \dd t.\]
	We can first eliminate $H$ be means of Lemma~\ref{lem:pt-H-estimate} via
	\[T(\cE(T) + \cD(T)) \lesssim (1+T) \!\int_0^T\!\! \int_\Omega \abs{\pt E}^2 \dd x \dd t
		+ \int_0^T\!\! \int_\Omega \sigma E \cdot E\dd x \dd t.\]
	Here we will insert the decomposition $\pt E = \pt V_h + \pt V_{i} + \pt \nabla p$ established before Lemma~\ref{lem:energy}.
	Next, the homogeneous system~\eqref{eq:V_h} for $(V_h,W_h)$, Lemma~\ref{lem:energy}, and
	Remark~\ref{rem:approx} show that the integrals
	\[\int_\Omega \big(\epsilon \pt V_h \cdot \pt V_h + \mu \pt W_h \cdot \pt W_h \big)\dd x
		= \int_\Omega \big(\epsilon \pt V_h \cdot \pt V_h + \inv\mu\curl V_h \cdot \curl V_h \big)\dd x \]
	are constant in time. As $T\ge 1$, we thus obtain
	\begin{align*}
		T(\cE(T) + \cD(T)) & \lesssim T \!\int_0^T\!\!\int_\Omega \big(|\pt V_h|^2 + |\pt V_{i}|^2 + |\pt \nabla p|^2\big)\dd x \dd t
		+ \int_0^T \!\!\int_\Omega \sigma E \cdot E\dd x \dd t \\
		 & \lesssim T^2\! \int_\Omega \big(|\pt V_h(0)|^2 +|\curl V_h(0)|^2 \big)\dd x \\
		 & \seq+ T \int_0^T\!\!\int_\Omega \big(|\pt \nabla p|^2 + |\pt V_{i}|^2 \big)\dd x \dd t
		+ \int_0^T\!\! \int_\Omega \sigma E \cdot E\dd x \dd t.
	\end{align*}
	The observability estimate of Proposition~\ref{prop:hom_obs} and Remark~\ref{rem:approx} yield
	\begin{align*}
		T(\cE(T) +\cD(T)) & \lesssim T^2\! \int_0^T \!\!\int_\omega \pt V_h \cdot \pt V_h \dd x \dd t
		+ T\! \int_0^T\!\!\int_\Omega \big(\abs{\pt \nabla p}^2 +\abs{\pt V_{i}}^2 \big)\dd x \dd t \\
		 & \seq + \int_0^T \int_\Omega \sigma E \cdot E\dd x \dd t
	\end{align*}
	for $T\geq T_O$. After replacing again $V_h = E - V_{i} - \nabla p$, condition \eqref{eq:conditions-on-conductivity} leads to
	\begin{align*}
		T(\cE(T) + \cD(T)) & \lesssim T^2 \!\!\int_0^T \!\!\! \int_\Omega \!\sigma \pt E(t) \cdot \pt E(t) \dd x \dd t
		+T^2\!\! \int_0^T \!\!\!\int_\Omega\!\! \big(|\pt \nabla p|^2 \!+ |\pt V_{i}|^2\big) \textrm{d} x \dd t \\
		 & \seq + \int_0^T \!\!\int_\Omega \sigma E \cdot E\dd x \dd t.
	\end{align*}
	Using Lemmas~\ref{lem:pt-V_nh-estimate}, \ref{lem:reg_of_ptp} and \ref{lem:energy-decreasing}, we finally deduce
	\begin{align*}
		T(\cE(T) + \cD(T)) & \lesssim T^4\! \int_0^T\!\!\int_\Omega \big(\sigma E\cdot E+\sigma \pt E \cdot \pt E\big)\dd x \dd t \\
		 & \leq C T^4 \big(\cE(0) - \cE(T) + \cD(0) - \cD(T)\big).
	\end{align*}
	Setting $\gamma=\frac{CT^4}{CT^4+T} <1$, we conclude
	\[\cE(T) + \cD(T) \leq \gamma \big(\cE(0) + \cD(0)\big). \qedhere \]
\end{proof}

Theorem~\ref{thm:exp-decay} now follows by a simple argument.

\begin{proof}[Proof of Theorem~\ref{thm:exp-decay}]
	First let $(E_0,H_0)\in D(A)$. Iterating the estimate from Proposition~\ref{prop:decay}, we obtain constants $\tilde M \ge 1$
	and $\omega>0$ such that
	\[ \cE(t) + \cD(t) \leq \tilde M\e^{-\omega t/2} \big(\cE(0) + \cD(0)\big)\ \]
	for all $t\ge 0$. Since $\cE(t) =\|T(t)(E_0,H_0)\|_X^2$ and $\cD(t) = \|AT(t)(E_0,H_0)\|_X^2$ by \eqref{def:ED},
	we have shown that $T(\cdot)$ exponentially decays in $D(A)$ with the graph norm. This space is isomorphic to $X$ by
	$(I-A)^{-1}$,	so that the assertion follows.
\end{proof}

Finally, we remove the divergence constraints in Theorem~\ref{thm:exp-decay} by projecting in
$X_e=L_{\epsilon}^2(\Omega) \times L^2_\mu(\Omega)$ onto $\Nu(A_e)^\perp=\ol{\Ra(A_e)}$
for the extension $A_e$ of $A$, see the discussion after \eqref{eq:solution}.
Here we proceed similar to \cite{NS25}. Then the theorem will
imply the exponential decay of the extended semigroup $T_e(\cdot)$ to the kernel $\Nu(A_e)$.

\begin{Lemma}\label{lem:tildeA}
	Let \ref{hypothesis} and \eqref{eq:conditions-on-conductivity} hold. We then have
	\[ \Nu(A_e) = \big\{(E,H)\in H_{t0}(\curl0)\times H(\curl 0)\,\big|\,E=0 \text{ on } \omega\big\}=\Ra(A_e)^\perp.\]
	Moreover, the orthogonal projection $P$ onto $\Nu( A_e)$ commutes with $T_e(\cdot)$, $T_e(t)P=I$ for $t\ge1$, and
	$\Nu( A_e)^\perp=\ol{\Ra(A_e)}$ is contained in $X$.
\end{Lemma}

\begin{proof}
	\begin{enumerate}[label = \arabic*), wide, labelwidth=0pt]
		\item \label{proof:tildeA-1} Take $w=(E,H)\in \Nu(A_e)$. This means that $\curl E=0$ and $\curl H=\sigma E$.
		 As $\tr_t E=0$, integration by parts yields
		 \[ 0=(A_ew|w)_{X_e}= \int_\Omega\big(\curl H\cdot E -\sigma E\cdot E -\curl E\cdot H\big)\dd x
			 = -\int_\Omega |\sigma^\frac12 E|^2\dd x, \]
		 so that $\sigma E=0$ and thus $E=0$ on $\omega$ by \eqref{eq:conditions-on-conductivity}. This shows `$\subseteq$'
		 in the first asserted identity. The converse inclusion is clear because of $\supp \sigma =\ol{\omega }.$

		\item Step~\ref{proof:tildeA-1} and integration by parts imply that the kernel $\Nu(A_e)$ is orthogonal to
		 the range $\Ra(A_e)$. To show $\Nu(A_e)= \Ra(A_e)^\perp$, take $h=(f,g)\in X_e$ with
		 \[0=(h|A_e w)_{X_e}= \int_\Omega\big(f\cdot \curl H - f\cdot \sigma E-g\cdot \curl E\big)\dd x\]
		 for all $w=(E,H)\in D(A_e)$. Choosing $(0,H)$ and $(E,0)$ with $E,H\in H^1_0(\Omega)$, we see that $\curl f=0$ and
		 $\sigma f+\curl g=0$. If we insert $(0,H)$ with $H\in H^1(\Omega)$, it follows
		 $\langle \tr_t f, H\rangle_{H^{-1/2}(\partial \Omega)}=0$ so that $h\in D(A_e)$. The formula in display with $h=w$
		 then implies $\int_\omega \sigma f\cdot f\dd x =0$ which yields $f=0$ on $\omega$ and $\curl g=0$; i.e.,
		 $h\in \Nu(A_e)$ by step~\ref{proof:tildeA-1} as needed.

		\item Hence, $\ol{\Ra(A_e)}=\Nu(A_e)^\perp$ is the kernel of $P$, implying $PA_e=0=A_eP$ on $D(A_e)$,
		 and thus $PT_e(t)=T_e(t)P=P$ for $t\ge0$. Let $(f,g)=A_e(E,H)$ for some $(E,H)\in D(A_e)$.
		 Then $g=-\mu ^{-1}\curl E$ belongs to $H^\Sigma_{n 0}(\div_\mu 0)$ by Proposition~6.1.4 in \cite{ciarlet18}, and
		 $\epsilon f_{|\upsilon} =\curl H$ to $H^\Gamma(\div 0,\upsilon)$ because of
		 \eqref{eq:conditions-on-conductivity}, Proposition~IX.1.3 in \cite{dautray1990}, and the density
		 of $H^1(\upsilon)$ in $H(\curl,\upsilon)$. As a result, $\ol{\Ra(A_e)}$ is contained in $X$. \qedhere
	\end{enumerate}
\end{proof}

\begin{Corollary}\label{cor:main}
	Let \ref{hypothesis}, \eqref{eq:conditions-on-conductivity}, and \eqref{eq:non-trapping} hold. Then there exist constants
	$M' \ge 1$ and $\omega>0$ such that for $(E_0,H_0) \in X_e=L_{\epsilon}^2(\Omega) \times L^2_\mu(\Omega)$ we have
	\[ \norm{(T_e(t) - P)(E_0,H_0)}_{L^2} \leq M'\e^{-\omega t} \norm{(E_0,H_0)}_{L^2}, \qquad t\geq 0.\]
\end{Corollary}
\begin{proof}
	Lemma~\ref{lem:tildeA} yields $T_e(t) - P = T_e(t)(I-P) = T(t)(I-P)$ so that
	the result follows from Theorem~\ref{thm:exp-decay}.
\end{proof}




\bibliographystyle{plain}
\bibliography{sources}

\end{document}